\documentclass[12pt]{amsart}
\usepackage{amsmath}
\usepackage{amssymb}
\usepackage{psfig}  

\catcode`@=11 \@mparswitchfalse  
\newcounter{mnotecount}[section]

\date{\today} 

\newcommand{\bi}{\bf\itshape}  

\numberwithin{equation}{section}   

\newcommand{\la}{\langle}	   	
\newcommand{\ra}{\rangle}          	

\renewcommand{\(}{\left(}
\renewcommand{\)}{\right)}

\newcommand{\R}{{\mathbf R}}		 
\newcommand{\e}{\varepsilon}             
\newcommand{\nothing}{\varnothing}       

\newcommand{\nn}{\nonumber}	   	
\newcommand{\note}[1]{}                    

\theoremstyle{plain}

\swapnumbers 

\newtheorem{thm}{Theorem}[section]
\newtheorem{lemma}[thm]{Lemma}
\newtheorem{prop}[thm]{Proposition}
\newtheorem{cor}[thm]{Corollary}

\theoremstyle{definition}
\newtheorem{defn}[thm]{Definition}
\newtheorem{adefn}[thm]{Alternative Definition}

\theoremstyle{remark}
\newtheorem{remark}[thm]{Remark}

\newtheorem{example}[thm]{Example}

\renewcommand{\phi}{\varphi}				
\newcommand{\cd}{,\dots,}	
\newcommand{\ol}{\overline}
\newcommand{\dist}{\operatorname{dist}\nolimits}
\newcommand{\cn}{\colon}	

\newcommand{\f}{\partial}
\newcommand{\fp}[1]{\boldsymbol{\{}#1\boldsymbol{\}} }	
\newcommand{\Leb}{{\mathcal L}}			
	
\newcommand{\bpi}{{\boldsymbol\pi}}		
\newcommand{\alg}{{\mathcal A}}			
\newcommand{\N}{{\mathbf N}}	

\newcommand{\prob}{{\mathcal P}}
\newcommand{\almost}{\operatorname{AlmCon}}
\renewcommand{\setminus}{\smallsetminus}
\newcommand{\eq}[1]{~\eqref{#1}}
\newcommand{\mean}{{\mathcal M}}
\newcommand{\nodes}{{\mathcal N}}

\newcommand{\sym}{\operatorname{sym}}

\title[Almost Convex Functions]{A General Theory of Almost Convex Functions.}

\author[Dilworth]{S. J. Dilworth}
\address{\null\hskip-1\parindent Department of Mathematics\newline
University of South Carolina\newline
Columbia, S.C. 29208, USA\newline 
\null\newline
{\tt dilworth\char'100math.sc.edu}\newline 
{\tt howard\char'100math.sc.edu}\newline
{\tt roberts\char'100math.sc.edu}
}

\author[Howard]{Ralph Howard}

\author[Roberts]{James W. Roberts}

\thanks{The research of the second author was supported in part from ONR Grant
N00014-90-J-1343 and ARPA-DEPSCoR Grant DAA04-96-1-0326}

\begin{document}
\allowdisplaybreaks 

\begin{abstract}
Let $\Delta_m=\{(t_0\cd t_m)\in \R^{m+1}: t_i\ge 0, \sum_{i=0}^mt_i=1\}$ be
the standard $m$-dimensional simplex.  Let $\nothing\ne S\subset
\bigcup_{m=1}^\infty\Delta_m$, then a function $h\cn C\to \R$ with
domain a convex set in a real vector space is {\bi $S$-almost convex\/}
iff for all $(t_0\cd t_m)\in S$ and  $x_0\cd x_m\in C$ the inequality
$$
h(t_0x_0+\cdots+t_mx_m)\le 1+ t_0h(x_0)+\cdots+t_mh(x_m)
$$
holds. A detailed study of the properties of $S$-almost convex
functions is made.  If $S$ contains at least one
point that is not a vertex, then an extremal $S$-almost convex
function $E_S\cn \Delta_n\to \R$ is constructed with the properties
that it vanishes on the vertices of $\Delta_n$ and if $h\cn
\Delta_n\to \R$ is any bounded $S$-almost convex function with
$h(e_k)\le 0$ on the vertices of $\Delta_n$, then $h(x)\le E_S(x)$ for
all $x\in \Delta_n$.  In the special case 
$S=\{(1/(m+1)\cd 1/(m+1))\}$, the barycenter of $\Delta_m$, very
explicit formulas are given for $E_S$ and
$\kappa_S(n)=\sup_{x\in\Delta_n}E_S(x)$.  These are of interest as
$E_S$ and $\kappa_S(n)$ are extremal in various geometric and analytic
inequalities and theorems. 
\end{abstract}

\maketitle

\tableofcontents

\section{Introduction.}

Let $C$ be a convex set in a real vector space and let $h\cn C\to R$.
Then according to Hyers and Ulam~\cite{Hyers-Ulam:convex} for $\e>0$,
$h$ is $\e$-approximately convex iff
\begin{equation}\label{hyers-ulam}
h((1-t)x+ty)\le \e+ (1-t)h(x)+th(y),\quad \text{for all} \quad t\in [0,1].
\end{equation}
In~\cite{Hyers-Ulam:convex} they show that if $h$ is
$\e$-approximately convex and $C\subseteq\R^n$ then there is a convex
function $g\cn C\to \R$ and a constant $C(n)$ only depending on the
dimension so that $|h(x)-g(x)|\le \frac12C(n)\e$.  In a previous paper
we show the sharp constant is
$$
C(n)=\lfloor\log_2n\rfloor+\frac{2(n+1-2^{\lfloor
\log_2n\rfloor})}{n+1}.
$$
(Here $\lfloor\cdot\rfloor$ is the floor, or greatest integer
function, and $\lceil\cdot \rceil$ is the ceiling function, that is
$\lceil x\rceil$ is the smallest integer greater than or equal to
$x$.)  In the present paper we generalize the notion of approximate
convexity and give the sharp constants in the corresponding Hyers-Ulam
type theorems.  This is done by finding the extremal approximately
convex function on the simplex that vanishes on the vertices.

Let us put the these problems in a somewhat larger setting.  First, by
replacing $h$ by $\e^{-1}h$ in~\eqref{hyers-ulam}, there is no loss of
generality in assuming that $\e=1$.  Then many natural notions of
generalized convexity are covered in the following definition.  Let
$\Delta_m=\{(t_0\cd t_m)\in \R^{m+1}: t_i\ge 0,
\sum_{i=0}^mt_i=1\}$ be the standard $m$-dimensional simplex.

\begin{defn}\label{def:S-almost}
Let $V$ a vector space over the reals and let $\nothing\ne C\subseteq
V$ be a convex set and let $\nothing \ne S \subseteq
\bigcup_{m=1}^\infty\Delta_m$.  Then a function $h\cn C\to \R$ is
{\bi $S$-almost convex\/} on $C$ iff for all $(t_0\cd t_m)\in S$ and
$x_0\cd x_m\in C$ the inequality
$$
h\(\sum_{i=0}^mt_ix_i\)\le 1+\sum_{i=0}^mt_ih(x_i)
$$
holds.  We denote by
$$
\almost_S(C):=\{h: \text{$h$ is $S$-almost convex on $C$}\}
$$
the set of almost convex functions $h\cn C\to \R$.\qed
\end{defn}

The case of $S=\Delta_1$ corresponds to the case studied by Hyers and
Ulam~\cite{Hyers-Ulam:convex} and others (cf.~the
book~\cite{Hyers-Isac-Rassias} for more information and references).
When $S=\{(1/2,1/2)\}$ the $S$-almost convex functions are just the
functions that satisfy
$$
h\(\frac{x+y}{2}\)\le 1+\frac{h(x)+h(y)}{2}.
$$
which are the approximately midpoint convex functions, (sometimes
called the approximately Jensen convex functions) which also have been
studied by several authors.

We give a general theory of $S$-almost convex functions. In
particular when $S$ has at least one point that is not a vertex
we construct (Definition~\ref{def:E} and Theorem~\ref{E-extr}) a bounded
$S$-almost convex function $E_S^{\Delta_n}\cn \Delta_n\to \R$ such
that if $h\cn \Delta_n\to \R$ is bounded, $S$-almost convex, and
$h(e_k)\le 0$ on the vertices of $\Delta_n$ then $h(x)\le
E^{\Delta_n}_S(x)$ for all $x\in \Delta_n$.  Then the number
$\kappa_S(n):=\sup_{x\in \Delta_n}E^{\Delta_n}_S(x)$ is the sharp
constant in stability theorems of Hypers-Ulam type and the function
$E^{\Delta_m}_S$ is the function that shows it is sharp (See
Theorem~\ref{stability}.) 

Probably the most natural choices, for $S$ are $S=\Delta_m$, a
simplex, and $S=\{(1/(m+1)\cd 1/(m+1))\}$, the barycenter of a
simplex.  In these cases we are able to give very explicit formulas
both for the extremal function $E_S^{\Delta_n}$ and for the constant
$\kappa_S(n)=\sup_{x\in \Delta_n}E^{\Delta_n}_S(x)$.  (For the case
$S=\Delta_m$ this was done in our earlier
paper~\cite{Dilworth-Howard-Roberts:3} where
$$
\kappa_{\Delta_m}(n)=\lfloor\log_{m+1}n\rfloor+\frac{\lceil
(m+1)\((n+1)-(m+1)^{\lfloor \log_{m+1}n\rfloor}\)/m\rceil}{n+1}.
$$
For the case of $S$ the barycenter of $\Delta_m$ see
Theorem~\ref{main-bary}, where the value is given as
\begin{equation}\label{kappa(m)}
\kappa_{\{(1/(m+1)\cd 1/(m+1))\}}(n)=\lfloor \log_{m+1}n\rfloor+1+
\frac{n}{m(m+1)^{\lfloor \log_{m+1}n\rfloor}}.
\end{equation}
(This differs from the notation of Theorem~\ref{main-bary} by the
substitution $B=m+1$.))  There is an interesting dichotomy in these
two cases.  When $S=\Delta_m$ then $E_S^{\Delta_n}$ is a concave
piecewise linear function that is continuous on the interior
$\Delta_n^\circ$ of $\Delta_n$ and the maximum occurs at the
barycenter of $\Delta_n$.  (See~\cite{Dilworth-Howard-Roberts:3}.)
However when $S=\{(1/(m+1)\cd 1/(m+1))\}$ is the barycenter of
$\Delta_m$ then $E_S^{\Delta_n}$ is discontinuous on a dense subset of
$\Delta_n$ and the graph of $E_S^{\Delta_n}$ is a fractal with a large
number of self similarities and the maximum does \emph{not} occur at
the barycenter of $\Delta_m$.  See Figure~\ref{E-graph}.  We also note
the somewhat surprising fact, that, as functions of $n$, both
$\kappa_{\Delta_m}(n)$ and $\kappa_{\{(1/(m+1)\cd 1/(m+1))\}}(n)$ have
the same order of growth, i.e. $\lfloor\log_{m+1}n\rfloor+O(1)$.

This paper is not completely self-contained.  Several of the results have
proofs that are very similar to the proofs in our earlier
paper~\cite{Dilworth-Howard-Roberts:1} and at several places we refer
the reader to~\cite{Dilworth-Howard-Roberts:1} for proofs.

\subsection{Definition and basic properties.}

Let $\Delta_m:=\{(t_0\cd t_m): \sum_{k=0}^mt_k=1, t_k\ge 0\}$ be the
standard $m$-dimensional simplex.  For the rest of this section we fix
a subset
$$
S\subset \bigcup_{m=1}^\infty \Delta_m.
$$
It follows easily from the definition  of $S$-almost convex that
$\almost_S(C)$ is a convex subset of the vector space of all functions
from $C$ to $\R$.

It is useful to make a  distinction between two cases:

\begin{defn}\label{def:cases}
If $S\subseteq
\bigcup_{m=1}^\infty\Delta_m$ then
\begin{enumerate}
\item If $S\nsubseteq \bigcup_{m=1}^N\Delta_m$ for any finite $N$ then
$S$ is of {\bi infinite type\/}.
\item If  $S\subseteq \bigcup_{m=1}^N\Delta_m$ for some $N$ then $S$ is
of {\bi finite type\/}.  If further $S\subseteq \Delta_m$ for some $m$
then $S$ is {\bi homogeneous\/}.\qed
\end{enumerate}
\end{defn}

\begin{remark}\label{S-compact}
If we assume that the union $ \bigcup_{m=1}^\infty\Delta_m$ is
disjoint and has the natural topology ($U\subseteq
\bigcup_{m=1}^\infty\Delta_m$ is open iff $U\cap \Delta_m$
is open in $\Delta_m$ for all $m$) then it is not hard to see that $S$
is of finite type if and only if it has compact closure in $
\bigcup_{m=1}^\infty\Delta_m$.\qed
\end{remark}

When considering $S$-almost convex functions there is no real
distinction between $S$ of finite type and $S$ homogeneous.

\begin{prop}\label{1Delta}
Let $S\subseteq \bigcup_{m=1}^N\Delta_m$. For $m\le N$ let
$\iota^m_N\cn \Delta_m\to \Delta_N$ be the inclusion 
$\iota^n_N(t_0\cd t_m)=(t_0\cd t_m,0\cd 0)$ and set
$S_m^*=\iota^m_N[S\cap \Delta_m]\subseteq \Delta_N$.  Let
$S^*=\bigcup_{m=1}^NS_m^*\subseteq \Delta_N$.  Then for any convex
subset $C$ of a real vector space  $\almost_{S^*}(C)=\almost_S(C)$.
\end{prop}

\begin{proof}
This is a more or less straightforward chase though the definition.
\end{proof}

The proof of the following is also straightforward and left to the
reader.

\begin{prop} \label{symmetric}
Let $S\subseteq \Delta_m$ and let 
$$S^*=\bigcup_{\rho\in
\sym(m+1)}\{(t_{\rho(0)},t_{\rho(1)}\cd t_{\rho(m)}): (t_0,t_1\cd
t_m)\in S\}
$$ 
where $\sym(m+1)$ is the group of all permutations of
$\{0,1\cd m\}$.  Then for any convex
subset $C$ of a real vector space $\almost_{S^*}(C)=\almost_S(C)$.
\end{prop}

The following is also trivial.

\begin{prop}\label{incl}
Let $S_1\subseteq S_2\subseteq \bigcup_{m=1}^\infty\Delta_m$.  Then 
for any convex subset $C$ of a real vector space
$\almost_{S_2}(C)\subseteq \almost_{S_1}(C)$.\qed
\end{prop}

The following can be used to reduce certain questions about $S$-almost
convex functions to the case where $S\subseteq \Delta_1$.

\begin{prop}\label{S-part}
Let $S\subseteq \bigcup_{m=1}^\infty\Delta_m$ and let $S_1$ be a
nonempty subset of $S\cap\Delta_m$ for some $m$.  Let $N_0\cd
N_k$ be a partition of the set $\{0,1\cd m\}$ into $k+1$ nonempty sets
and let
$$
S_2:=\{ (\alpha_0(t),\alpha_2(t)\cd \alpha_k(t)): t\in S_1\}\subseteq \Delta_k
$$
where
$$
\alpha_j(t):=\sum_{i\in N_j}t_i.
$$
Then 
$$
 \almost_{S}(C)\subseteq\almost_{S_2}(C)
$$ 
for any convex subset $C$ of a real vector space.  In particular if
$(t_0\cd t_m)\in S$ and for some $k\in \{0\cd m-1\}$ we set
$\alpha=t_0+\cdots +t_k$ and $\beta=t_{k+1}+\cdots+t_m$ then any $S$
almost convex function $h$ will satisfy $h(\alpha x_0)+h(\beta x_1)\le
1+\alpha h(x_0)+\beta h(x_1)$.
\end{prop}

\begin{proof}
Let $C$ be a convex subset of a real vector space and let $y_0\cd
y_k\in C$,  $\alpha\in S_2$ and $ h\in  \almost_{S}(C)$.  Let
$x_0\cd x_m \in C$ be defined by 
$$
x_i=y_j \quad \text{if} \quad i\in N_j
$$
As $\alpha\in S_2$ there is a $t=(t_0\cd t_m)\in S_1\subseteq S$ so
that  $\alpha_j=\sum_{i\in N_j}t_i$.  Then as $ h$ is $S$-almost
convex
$$
 h\bigg(\sum_{j=0}^k\alpha_jy_j\bigg)=
 h\bigg(\sum_{i=0}^mt_ix_i\bigg)\le 1+ \sum_{i=0}^kt_i h(x_i)
=1+\sum_{i=0}^m\alpha_j h(y_j).
$$
Thus $ h\in \almost_{S_2}(C)$.
\end{proof}

It is useful to understand when an $S$-almost convex function is bounded. 

\begin{thm}\label{leb-bded}
Let $S\subseteq \bigcup_{m=1}^\infty\Delta_m$ and assume that $S$
contains at least one point that is not a vertex (that is there is
$(t_0\cd t_m)\in S$ with $\max_{i}t_i<1$). Let $U$ be a convex open
set in $\R^n$.  Then any $S$-almost convex function $ h\cn U\to \R$
which is Lebesgue measurable is bounded above and below on any compact
subset of $U$.
\end{thm}

\begin{proof}
Let $(t_0\cd t_m)\in S$ with $\max_{i}t_i<1$.  Then there is a $k\in
\{0\cd m-1\}$
so that if $\alpha=t_0+\cdots+t_k$ and
$\beta=t_{k+1}+\cdots+t_m$, then $0<\alpha,\beta<1$, 
$\alpha+\beta=1$ and by Proposition~\ref{S-part}
$$
 h(\alpha x_0+\beta x_1)\le 1+\alpha h(x_0)+\beta h(x_1).
$$
We assume that $\alpha\le \beta$, the case of $\alpha>\beta$ having a
similar proof.  As any compact subset of $U$ is contained in a bounded
convex open subset of $U$ we can also assume, without loss of
generality, that $U$ is bounded.

Let $K\subset U$ be compact and let $r=\dist(K,\f U)$.  For any $x\in
\R^n$ let $B_r(x)$ be the open ball of radius $r$ about $x$.  Then for
any $a\in K$ we have $B_r(a)\subseteq U$.  For $a\in K$ define
$\theta_a\cn\R^n\to \R^n$ by
$$
\theta_a(x)=\frac1\beta a-\frac\alpha\beta x.
$$
Then it is easy to check that $\theta_a(a)=a$ for all $a\in \R^n$ and
$\alpha x+\beta\theta_a(x)=a$ for all $x\in \R^n$.  Also $\theta_a$ is
a dilation in the sense that
$\|\theta_a(x_1)-\theta_a(x_0)\|=(\alpha/\beta)\|x_1-x_0\|$ for all
$x_0,x_1\in \R^n$.  As $\theta_a(a)=a$ and $(\alpha/\beta)\le 1$ this
implies $\theta_a[B_a(r)]=B_a((\alpha/\beta)r)\subseteq B_a(r)$.  Let
$\Leb^n$ be Lebesgue measure on $\R^n$.  Then for any measurable subset $P$ of
$\R^n$
$$
\Leb^n(\theta_a[P])=(\alpha/\beta)^n\Leb^n(P).
$$
Choose a positive real number $\e$ so that
\begin{equation}\label{e-def}
\(1+\(\frac\alpha\beta\)^n\)\e < \(\frac\alpha\beta\)^n\Leb^n(B(r))
\end{equation}
where $B(r)$ is the open ball of radius $r$ about the origin.  Because
$h$ is measurable and $\Leb^n(U)<\infty$ there is a positive $M$ so
large that
$$
\Leb^n\{x\in U:  h(x)> M\} < \e.
$$
Therefore if $V:=\{x\in U:  h(x)\le M\}$ then $\Leb^n(U\setminus
V)<\e$.  Let $A:=B_a(r)\cap V$.  We now claim that $A\cap \theta_a[A]$
has positive measure.  For if not then $A$ and $\theta_a[A]$ would be
essentially disjoint subsets of $B_r(a)$ and therefore, using that
$\Leb^n(\theta_a[A])=(\alpha/\beta)^n\Leb^n(A)$,
\begin{align*}
\Leb^n(B_a(r))&\ge
\Leb^n(A)+\Leb^n(\theta_a[A])\\
&=\(1+\(\frac\alpha\beta\)^n\)\Leb^n(A)\\
&\ge \(1+\(\frac\alpha\beta\)^n\)(\Leb^n(B_a(r))-\e)
\end{align*}
which can be rearranged as 
$\(1+(\alpha/\beta)^n\)\e \ge (\alpha/\beta)^n\Leb^n(B(r))$
contradicting~\eqref{e-def}.  Therefore $\Leb^n(A\cap
\theta_a[A])>0$ as claimed.  Let $a\ne x\in A\cap
\theta_a[A])$.  Then $x$ and $\theta_a(x)$ are both in $A=B_a(r)\cap
V$ and therefore $ h(x), h(\theta_a(x))\le M$.  Thus 
$$
 h(a)= h(\alpha x+\beta\theta_a(x))\le
1+\alpha h(x)+\beta h(\theta_a(x))\le 1+\alpha M+\beta M=M+1
$$
which shows that $h$ is bounded above on $K$.

To show that $h$ has a lower bound on compact subsets of $U$, let
$a\in U$ and let $r>0$ be small enough that the closed ball
$\ol{B}_a(r)$ is contained in $U$.  Then $\ol{B}_a(r)$is compact so by
what we have just done there is a constant $C>0$ so that $ h(x)\le C$ for
all $x\in B_a(r)$.  Let $x\in B_a(r)$.  Then, again as above,
$\theta_a(x)\in B_a(r)$, and therefore
$$
 h(a)=h(\alpha x+\beta\theta_a(x))\le
1+\alpha h(x)+\beta h(\theta_a(x)) \le 1+\alpha h(x)+\beta C
$$
which can be solved for $ h(x)$ to give
$$
 h(x)\ge \frac1\alpha( h(a)-1-\beta C).
$$
Therefore $ h$ is bounded below on $B_a(r)$.  But any compact subset
of $U$ can be covered by a finite number of such open balls and thus
$ h$ is bounded below on all compact subsets of $U$.
\end{proof}

The following will be needed later.

\begin{cor}\label{1D-bded}
Let $h\cn [a,b]\to \R$ be a Lebesgue measurable function so that
$ h(\alpha x+\beta y)\le 1+\alpha h(x)+\beta h(y)$ for some
$\alpha,\beta>0$ with $\alpha+\beta=1$ (that is $ h$ is $S$-almost
convex with $S=\{(\alpha,\beta)\}\subset\Delta_1$).  Then $ h$ is
bounded above on $[a,b]$.
\end{cor}

\begin{proof}
By doing a linear change of variable (which preserves $S$-almost
convexity) we can assume that $[a,b]=[0,1]$.  Also by replacing $ h$
by $x\mapsto h(x)-((1-x) h(0)+x h(1))$ we can assume that $ h(0)=
h(1)=0$.  Let $\delta=\alpha/(1+\alpha)$.  Then by
Theorem~\ref{leb-bded} there is a constant $C_1>0$ such that $ h(x)\le
C_1$ on $[\delta,1-\delta]$.  Let
$$
C_2=\max\{C_1,1/(1-\alpha)+\alpha C_1\}.
$$
We now show that $ h\le C_2$ on $[0,1]$.  If $x=0$, $x=1$, or $x\in
[\delta,1-\delta]$ this is clear.  Let $x\in (0,\delta)$ then the
choice of $\delta$ ensures that there is a $y\in [\delta,1-\delta]$ such that
$x=\alpha^ky$ for some positive integer $k$.  Also, as $y\in
[\delta,1-\delta]$,  $ h(y)\le C_1$. Therefore
\begin{align*}
 h(x)&= h(\alpha^ky)= h(\beta 0+\alpha \alpha^{k-1}y)\\
&\le 1+ \beta h(0)
	+\alpha h(\alpha^{k-1}y)=1+ \alpha h(\alpha^{k-1}y)\\
&\le 1+\alpha\big(1+\alpha h(\alpha^{k-2}y)\big)
	= 1+\alpha +\alpha^2 h(\alpha^{k-2}y)\\
&\le 1+\alpha+\alpha^2+\cdots+\alpha^{k-1}+\alpha^k h(y)\\
&\le \frac{1}{1-\alpha}+\alpha C_1\le C_2.
\end{align*}
If $x\in (1-\delta,1)$ a similar calculation shows that $ h(x)\le
C_2$ (or this can be reduced to the case $x\in (0,\delta)$ by the
change of variable $x\mapsto (1-x)$).  This completes the proof.
\end{proof}

\subsection{A general construction for the extremal $S$ almost convex
function on a simplex.}\label{sec:gen-extr}

We will show that on the $n$-dimensional simplex $\Delta_n$ there is a
pointwise largest bounded $S$-almost convex function that vanishes on
the vertices of $\Delta_m$.  We start with some definitions.

\begin{figure}[ht]
\centering
\begin{picture}(260,150)
\put(80,150){\circle*{4}}
\put(40,100){\circle*{4}}
\put(120,100){\circle*{4}}
\put(20,50){\circle*{4}}
\put(60,50){\circle*{4}}
\put(90,50){\circle*{4}}
\put(120,50){\circle*{4}}
\put(150,50){\circle*{4}}
\thicklines
\put(40,100){\line(4,5){40}}
\put(120,100){\line(-4,5){40}}
\put(20,50){\line(2,5){20}}
\put(60,50){\line(-2,5){20}}
\put(90,50){\line(3,5){30}}
\put(150,50){\line(-3,5){30}}
\put(120,50){\line(0,1){50}}

\bezier{10}(20,50)(10,25)(0,0)
\bezier{10}(20,50)(20,25)(20,0)
\bezier{10}(20,50)(30,25)(40,0)

\bezier{10}(60,50)(55,25)(50,0)
\bezier{10}(60,50)(65,25)(70,0)

\bezier{10}(90,50)(85,25)(80,0)
\bezier{10}(90,50)(95,25)(100,0)

\bezier{10}(120,50)(115,25)(110,0)
\bezier{10}(120,50)(125,25)(130,0)
\bezier{10}(120,50)(120,25)(120,0)

\bezier{10}(150,50)(145,25)(140,0)
\bezier{10}(150,50)(155,25)(160,0)
\bezier{10}(150,50)(150,25)(150,0)


\put(48,125){\footnotesize$t_0$}
\put(103,125){\footnotesize$t_1$}

\put(17,70){\footnotesize$s_0$}
\put(54,70){\footnotesize$s_1$}

\put(90,70){\footnotesize$r_0$}
\put(142,70){\footnotesize$r_1$}
\put(122,70){\footnotesize$r_2$}


\put(100,147){\footnotesize The root $:=$ unique node of rank $0$.}
\put(140,97){\footnotesize The rank one nodes.}
\put(170,47){\footnotesize The rank two nodes.}

\end{picture}

\caption[]{\footnotesize 
An $S$ ranked tree showing the labeling of the edges out of the root
by $t=(t_0,t_1)\in S$ and the edges out of the rank one nodes by
$s=(s_0,s_1)\in S$ and $r=(r_0,r_1,r_2)\in S$.  In our definition each
node will have at least two edges leaving it and the sum of the
weights $t_0\cd t_m$ of the weights of all edges leaving a node is
unity (as $(t_0\cd t_m)\in \Delta_m$).  Finally, in the definition of
tree used here, all branches are of infinite length.}
\label{fig:S-tree}
\end{figure}
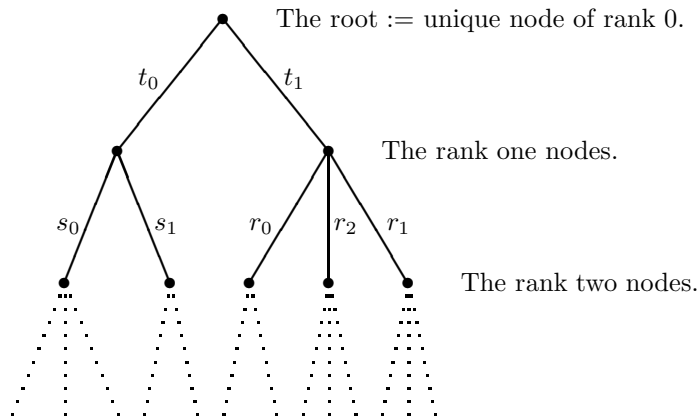

\begin{defn}\label{tree}
A {\bi tree\/}, $T$, is a collection of points $\nodes$, called {\bi
nodes\/}, and a set of (directed) {\bi edges\/} connecting some pairs
of nodes with the following properties: The set $\nodes$ is a disjoint
union $\nodes=\bigcup_{k=0}^\infty\nodes_k$ where $\nodes_0$ contains
exactly one point, the {\bi root of the tree}, each $\nodes_k$ is a
finite set and if $\nodes_k=\{v_1\cd v_m\}$ then $\nodes_{m+1}$ is a
disjoint union $\nodes_{m+1}=\mathcal{P}_1\cup \dots \cup
\mathcal{P}_m$ of nonempty sets where $\mathcal{P}_i$ is the set of
{\bi successors of $v_i$}.  The (directed) edges of the tree leave a
node and connect it to its successors and there are no other edges in
the tree (cf.~Figure~\ref{fig:S-tree}).  If $v$ is a node of the tree
then $r(v):=k$ where $v\in \nodes_k$ is the {\bi rank\/} of $v$.  A
{\bi branch\/} of the tree is a sequence of nodes $\la
v_k\ra_{k=0}^\infty$ where $v_0$ is the root, $r(v_k)=k$, and there is
an edge from $v_k$ to $v_{k+1}$.\qed
\end{defn}

We now consider trees with  extra structure, a labeling of the
edges in a way that will be used in defining the extremal
$S$-almost convex function.

\begin{defn}\label{def:ranked-tree}
Let $S\subseteq \bigcup_{m=1}^\infty \Delta_m$ be nonempty.  Then an
{\bi $S$-ranked tree\/} is a tree $T$ with its edges labeled by
non-negative real numbers in such a way that for any node $v$ of the
tree there is an element $t=(t_0\cd t_m)\in S$ so that there are
exactly $m+1$ edges leaving $v$ and these are labeled by $t_0\cd t_m$.
The number $t_i$ is the {\bi weight\/} of the edge it labels.  
Figure~\ref{fig:S-tree} shows
a typical $S$-ranked tree.\qed
\end{defn}

We now describe how an $S$-ranked tree determines a probability
measure on the set of branches of the tree.  Let $T$ be an $S$-ranked
tree and let $X=X(T)$ be the set of all branches of $T$.  If $\la
v_k\ra_{k=0}^\infty, \la w_k\ra_{k=0}^\infty\in X$ are two elements of
$X$ we can define a distance between them as $d(\la v_k\ra_{k=0}^\infty,
\la w_k\ra_{k=0}^\infty)=2^{-\ell}$ where $\ell$ is the smallest index
with $v_\ell\ne w_\ell$ (and $d(\la v_k\ra_{k=0}^\infty,
\la w_k\ra_{k=0}^\infty)=0$ if $\la v_k\ra_{k=0}^\infty=
\la w_k\ra_{k=0}^\infty$).  While we will not need to use this fact, it is
not hard to check that this makes $X$ into a compact metric space
which is homeomorphic to the Cantor set.

\begin{defn}\label{def:S-ranked-meas}
Let $S\subseteq \bigcup_{m=1}^\infty\Delta_m$ be nonempty and let $T$
be an $S$-ranked tree.  Then $T$ defines a measure on $X$, the set of
branches of $T$, as follows.  For $v$ a node of $T$ let $I(v)$ be the
set of branches of $T$ that pass through $v$.  If $k=r(v)$ is the
rank of $v$ then let $\la v_0,v_1\cd v_k\ra$ be the initial segment of
a branch passing through $v$ (so that $v=v_k$) and for $1\le i\le k$
let $s_i$ be the weight of the edge from $v_{i-1}$ to $v_i$.  Then
$\mu$ is the measure on $X$ such that
$$
\mu(I(v))=s_0s_1\cdots s_k.
$$
(That is $\mu(I(v))$ is the product of the weights of the edges along
an initial segment of a branch connecting the root to $v$.)  A measure
arising in this way will be called an {\bi $S$-ranked probability
measure\/}.\qed
\end{defn}

It follows from this definition that if $v$ is a node of $T$ and
$v_0\cd v_m$ are the successors of $v$ and $t=(t_0\cd t_m)\in S$
labels the edges from $v$ in such a way that $t_i$ labels the edge
from $v$ to $v_i$ then
$$
\mu(I(v_i))=t_i\mu(I(v)).
$$

It is useful to give a description of an $S$-ranked probability
measure that does not rely directly on its construction from an
$S$-ranked tree.  

\begin{adefn}\label{def:S-meas}
An {\bi $S$-ranked probability measure\/} is an ordered triple
$(X,\mu,\bpi)$ where $X$ is a nonempty set, $\bpi=\la
\pi_0,\pi_1,\pi_2,\dots\ra$ a sequence of finite partitions of $X$
into nonempty subsets such that $\pi_0=\{X\}$ and $\pi_{k+1}$ refines
$\pi_k$, $\mu$ is a measure defined on the $\sigma$-algebra,
$\alg(\bpi)$, generated by $\bigcup_{k=0}^\infty\pi_k$ so that for all
$j\ge 0$ and all $I\in \pi_j$, there exists $(t_0\cd t_m)\in S$ such
that if
$$
\{J\in \pi_{j+1}: J\subset I\} =\{I_0, I_1\cd I_m\}
$$
then
$$
\mu(I_i)=t_i\mu(I), \quad \text{for}\quad 0\le i\le m.
$$
If $I\in \bigcup_{k=0}^\infty\pi_k$ then the {\bi rank\/} of $I$ is
$r(I)=k$ where $I\in \pi_k$.  (The union $\bigcup_{k=0}^\infty\pi_k$
is disjoint so this is well defined.)\qed
\end{adefn}

Given an $S$-ranked probability measure $(X,\mu,\bpi)$ we can
construct an $S$-ranked tree by using for the set of nodes of the tree
$\nodes=\bigcup_{j=0}^\infty\pi_j$, letting $\nodes_k=\pi_k$ be the
set of nodes of rank $k$.  There is an edge from $I\in
\nodes_{j}=\pi_j$ to $J\in \nodes_{j+1}$ iff $J\subset I$ in this case
the weight of this edge is the $t_i$ such that $\mu(J)=t_i\mu(I)$.  In
most of what follows we will work with the alternative definition of
$S$-ranked probability~\ref{def:S-meas}, but will think of any such
measure as being constructed from an $S$-ranked tree as above.

\begin{example}\label{example1}
Suppose $S$ consists of a single point $(t_0\cd t_m)$ in the interior
of $\Delta_m$ (so that each $t_i$ is positive).  Then there is  only
one $S$-ranked probability measure i.e.\ $\mu=$ the product measure on
$[m]^\N$ where $[m]=\{0,1\cd m\}$ and $\mu=\nu\times\nu\times\cdots$
and $\nu$ is given on $[m]$ by $\nu(\{i\})=t_i$.  This uniqueness is
clear when viewed in terms of $S$-ranked trees as when $S$ is a one
point set there is clearly only one $S$-ranked tree.\qed
\end{example}

\begin{remark}\label{note1}
Let $(t_0\cd t_m)\in S$ and for each $i$ with $0\le i\le m$, let
$(X_i,\mu^{(i)},\bpi^{(i)})$ be an $S$-ranked probability measure on a
set $X_i$ where we assume $X_i\cap X_j=\nothing$ for $i\ne j$. We let
$X=\coprod_{i=0}^mX_i$ (the disjoint union of the $X_i$) and let
$\pi_0=\{X\}$.  For $j\ge 1$, set $\pi_j=\bigcup_{i=1}^m
\pi^{(i)}_{j-1}$. (This gives $\pi_1:=\{X_0\cd X_m\}$.)  Define a
measure $\mu$ on $\alg(\bpi)$ by
$\mu(A)=\sum_{i=1}^mt_i\mu^{(i)}(A\cap X_i)$.  Then $(X,\mu,\bpi)$ is an
$S$-ranked probability measure. Note that if $I\in \pi^{(i)}_j$ then
$r_{\mu^{(i)}}(I)=j$ and $r_{\mu}(I)=j+1$.\qed
\end{remark}

\begin{defn}\label{def:divides}
If $x=(x_0\cd x_n)\in \Delta_n$ and $\alpha=\la
\alpha_i\ra_{i=1}^\infty$ is a probability sequence in $\ell_1^+$
(that is $\sum_{i=1}^\infty\alpha_i=1$ and $\alpha_i\ge 0$) then {\bi
$x$ divides $\alpha$}, written as $x\mid \alpha$, iff
$\N=\{1,2,\dots\}$ can be partitioned into sets $N_0,N_1\cd N_n$ such
that
$$
x_k=\sum_{i\in N_k}\alpha_i\quad \text{for} \quad k=0,1\cd n.
$$
\null\vskip-22pt\qed
\end{defn}

\begin{defn}\label{def:E}
Define $E=E_S^{\Delta_n}\cn \Delta_n\to \R$ by
$$
E(x)=\inf \sum_{j=1}^\infty r_\mu(I_j)\mu(I_j)
$$
where the infimum is taken over all $S$-ranked probability measures
$(X,\mu,\bpi)$ and all disjoint sequences $\la
I_j\ra_{j=1}^\infty\subset\{\nothing\}\cup\bigcup_{k=0}^\infty\pi_k$ with
\begin{equation}\label{Ij-cond}
\quad \sum_{j=1}^\infty
\mu(I_j)=1\quad\text{and}\quad x\mid \la \mu(I_j)\ra_{j=1}^\infty.
\end{equation}
(This can be rephrased using disjoint sequences $\la I_j\ra\subset
\bigcup_{k=0}^\infty\pi_k$ which are either finite or countable.  But
it is notationally more convenient to take a finite sequence 
$\la I_j\ra_{j=1}^m$ and extend it to a sequence $\la
I_j\ra_{j=1}^\infty$ with $I_j=\nothing$ for $j\ge m+1$.)\qed
\end{defn}

In much of what follows it will be clear that the domain of $E$ is
$\Delta_n$ and we will just write $E_S$ or just $E$ rather than
$E_S^{\Delta_n}$.

\begin{remark}\label{note2}
For each $S$-ranked probability measure $(X,\mu,\bpi)$ we let
$\alg_i$ denote the finite algebra with elements of $\pi_i$ as its
atoms.  Then in the last definition let $\la I_j\ra_{j=1}^\infty\subset
\bigcup_{k=0}^\infty\alg_k$ be a disjoint sequence so that\eq{Ij-cond}
holds and let $\N=N_0\cd N_n$ be a partition of $\N$ so that
$x_k=\sum_{j\in N_k}\mu(I_j)$.  Then set $A_{i\,k}=\cup\{I_j:
r(I_j)=i, j\in N_k\}$.  Then
$$
 \sum_{j=1}^\infty r_\mu(I_j)\mu(I_j)=\sum_{k=0}^n\sum_{i=0}^\infty
 i\mu(A_{i\,k}).
$$
Therefore we could also define $E(x)$ by
$$
E(x)=\inf \sum_{k=0}^n\sum_{i=0}^\infty
 i\mu(A_{i\,k})
$$
where the infimum is taken over all $S$-ranked probability measures,
and all disjoint sequences $\la A_{i\,k}\ra_{0\le k\le n,\, 0\le i}$ so
that 
$$
A_{i\,k}\in \alg_i\quad \text{and}\quad \sum_{i}\mu(A_{i\,k})=x_k.
$$
\null\vskip-32pt\qed
\end{remark}

The following sum will be used later in this section and in
Section~3.
The proof is left
to the reader.

\begin{lemma}\label{lemma:sum}
Let $a,x\in \R$ with $|x|<1$ and $k$ an integer. Then
\begin{align*}
\sum_{j=0}^\infty(a+j)x^{k+j}&=ax^k+(a+1)x^{k+1}+(a+2)x^{k+2}+\cdots\\
&=\frac{ax^k}{1-x}+\frac{x^{k+1}}{(1-x)^2}=\frac{ax^k+(1-a)x^{k+1}}{(1-x)^2}.
\end{align*}
\null\vskip-20pt\qed
\end{lemma}

\begin{prop}\label{prop1}
For any nonempty $S\subset \bigcup_{m=1}^\infty\Delta_m$ we have $E_S(e_k)=0$
for all vertices of $\Delta_n$ and if $x\in \Delta_n$ is not a vertex
then $E_S(x)\ge 1$.  If $S$ contains a point $(t_0\cd t_m)$ which is
not a vertex, i.e.\ $\e:=\max_{i}t_i<1$, then  $E_S$ is
bounded on $\Delta_n$ and in fact has the upper bound
$$
E_S(x)\le 1+\frac{(2\e-\e^2)(n+1)}{(1-\e)^2}
$$
on $\Delta_n$.  Thus if $\inf_{t\in S}\max_{i}t_i=0$ (for example when
$S=\bigcup_{m=1}^\infty\Delta_m$) then $E$ is given by $E(e_k)=0$ and $E(x)=1$
for $x\in \Delta_n$ and $x$ not a vertex.
\end{prop}

\begin{proof}
If $x$ is a vertex of $\Delta_m$, which without lost of generality
we can take to be $x=e_0$, then let $(X,\mu,\bpi)$ be any
$S$-ranked probability measure and let $I_1=X$ and $I_j=\nothing$ for
$j\ge 2$.  Partition $\N$ as $N_0=\{1\}$ and $N_1\cd N_n$ an
arbitrary partition of $\N\setminus\{0\}$.  Then $r(I_1)=r(X)=0$ and
$\mu(I_j)=\mu(\nothing) =0$ for $j\ge 2$ and therefore
$$
0\le E(e_0)\le \sum_{j=0}^\infty r(I_j)\mu(I_j)=0.
$$
Thus $E(e_0)=0$.

Now assume that $x$ is not a vertex and let $(X,\mu,\bpi)$ be an
$S$-ranked probability measure and $\la I_j \ra_{j=1}^\infty$ with
$\sum_{j=1}^\infty \mu(I_j)=1$ and $x\mid \la
\mu(I_j)\ra_{j=1}^\infty$.  Then as $x$ is not a vertex we have that
$x_k<1$ for $0\le k\le n$ and therefore $\mu(I_j)\le x_k<1$.  Thus
$I_j\ne X$ and therefore $r(I_j)\ge 1$.  This gives
$$
\sum_{j=1}^\infty r_\mu(I_j)\mu(I_j)\ge \sum_{j=1}^\infty \mu(I_j)=1.
$$
Taking an infimum then gives that $E(x)\ge 1$.

Now assume that $S$ contains a point that is not a vertex and note that if
$S_1\subset S_2$ then $E_{S_2}(x)\le E_{S_1}(x)$ for all $x$.  Thus it
suffices to show that $E_S(x)$ is bounded when $S$ is a single point
$(t_0\cd t_m)$ with $\e=\max_it_i<1$.  Suppose $(x_0,x_1\cd x_n)\in
\Delta_n$.  We let $\mu$ be the product measure as in
Example~\ref{example1} and we let $\alg_i:=\alg(\pi_i)$ as in
Remark~\ref{note2} and use the alternative definition of $E_S$ given
in Remark~\ref{note2}.  For each $k$, $0\le k\le n$, we select
inductively a set $A_{i\,k}\in \alg_i$ with $\la A_{i\,k}\ra_{i,k}$
pairwise disjoint such that
$$
x_k-\e^i\le \sum_{j=0}^i \mu(A_{j\,k})\le x_k.
$$
Note that if $I\in \pi_i$, then
$\mu(I)\le \e^i$.  We carry out the the inductive selection as
follows:  Let 
$$
\mathcal{I}_i:=\bigg\{I\in \pi_i: I\cap \bigg(
\bigcup_{k=0}^n\bigcup_{j=0}^{i-1}A_{j\,k}\bigg)=\nothing\bigg\}
=\{I_1,I_2\cd I_M\}.
$$
Then
\begin{equation}\label{sum=1}
1=\sum_{k=0}^n\sum_{j=0}^{i-1}\mu(A_{j\,k})+\sum_{s=1}^M\mu(I_s).
\end{equation}
If $\sum_{j=0}^{i-1}\mu(A_{j\,0})\ge x_0-\e^i$, let
$A_{i\,0}=\nothing$.  If $\sum_{j=0}^{i-1}\mu(A_{j\,0})<x_0-\e^i$ let
$s_0$ be the first integer such that 
$$
\sum_{j=0}^{i-1}\mu(A_{j\,0})+ \sum_{s=0}^{s_0}\mu(I_s)\ge x_0-\e^i.
$$
Since $\mu(I_{s_0})\le \e^i$,
$$
\sum_{j=0}^{i-1}\mu(A_{j\,0})+\sum_{s=0}^{s_0}\mu(I_s)\le x_0.
$$
Let $A_{i\,0}=\bigcup_{s=0}^{s_0}I_s$.  Continue choosing from
$\{I_{s_0+1}\cd I_M\}$ to obtain $A_{i\,1}\cd A_{i\,n}$.  Note that
by~\eqref{sum=1}, the supply of atoms in $\mathcal{I}_i$ is sufficient
to choose the sets $A_{i\,0},A_{i\,1}\cd A_{i\,n}$.  For $i\ge 2$ we
have
$$
x_k-\e^{i-1}\le
\sum_{j=0}^{i-1}\mu(A_{j\,k})\le\sum_{j=0}^{i}\mu(A_{j\,k}) \le x_k
$$
which implies $\mu(A_{i\,k})\le \e^{i-1}$ for $i\ge 2$.  As
$\mu(A_{1\,k})\le x_k$ we can use Lemma~\ref{lemma:sum} (with $a=0$)
to compute
$$
\sum_{i=0}^\infty i\mu(A_{i\,k})
=\mu(A_{1\,k})+\sum_{i=2}^\infty i\mu(A_{i\,k})
\le x_k+\sum_{i=2}^\infty i\e^{i-1}
=x_k+\frac{2\e-\e^2}{(1-\e)^2}.
$$
Thus, in the notation of Remark~\ref{note2},
$$
E(x)\le \sum_{k=0}^n\sum_{i=1}^\infty i\mu(A_{i\,k})\le 
\sum_{k=0}^nx_k+\frac{(n+1)\e}{(1-\e)^2}=1+\frac{(n+1)(2\e-\e^2)}{(1-\e)^2}
$$
which bounds $E$ as required.
\end{proof}

\begin{prop}\label{prop2}
The function $E=E_S$ is $S$-almost convex on $\Delta_n$.
\end{prop}

\begin{proof}
Let $(t_0,t_1\cd t_m)\in S$ and $y_0,y_1\cd y_m\in \Delta_n$.  For
$0\le i\le m$, let $(X_i, \mu^{(i)},\bpi^{(i)})$ be an $S$-ranked
probability measure.  We let $\la
I^{(i)}_j\ra_{j=1}^\infty\subset\{\nothing\}\cup\bigcup_{k=1}^\infty
\pi_k^{(i)}$ be a disjoint sequence such that $y_i\mid \la
\mu^{(i)}(I^{(i)}_j)\ra_{j=1}^\infty$.  Now let $\mu$ be the
$S$-ranked probability measure on $X=\coprod_{i=0}^mX_i$ as in
Remark~\ref{note1}, i.e.\ $\mu(A)=\sum_{i=0}^mt_i\mu^{(i)}(X_i\cap
A)$.  It is easily checked that 
$\sum_{i=0}^mt_iy_i\mid \la \mu(I^{(i)}_j)\ra_{i=1}^m{}_{j=1}^{\infty}$ (and
$\sum_{i,j}\mu(I^{(i)}_j)=1$).  Thus
\begin{align*}
E\bigg(\sum_{i=0}^mt_iy_i\bigg)
&\le \sum_{i=0}^m\sum_{j=1}^\infty r_\mu\big(I^{(i)}_j\big)\mu\big(I^{(i)}_j\big)\\
&=\sum_{i=0}^m\sum_{j=1}^\infty
r_\mu\big(I^{(i)}_j\big)t_i\mu^{(i)}\big(I^{(i)}_j\big)\\
&=\sum_{i=0}^m\sum_{j=1}^\infty
\big[r_{\mu^{(i)}}\big(I^{(i)}_j\big)+1\big]t_i\mu^{(i)}\big(I^{(i)}_j\big)\\
&=1+\sum_{i=0}^m t_i\bigg(\sum_{j=1}^\infty
\big(r_{\mu^{(i)}}I^{(i)}_j\big)\mu^{(i)}\big(I^{(i)}_j\big)\bigg).
\end{align*}
Taking the infimum over all $\mu^{(0)}\cd \mu^{(m)}$ on the right hand side of
this gives $E\(\sum_{i=0}^mt_iy_i\)\le 1+ \sum_{i=0}^mE(y_i)$
which completes the proof.
\end{proof}

\begin{thm}\label{E-extr}
The function $E=E_S$ is the extremal $S$-almost convex function on
$\Delta_n$ in the sense that if $ h$ is a bounded $S$-almost convex
function on $\Delta_n$ with $ h(e_k)\le 0$ for $0\le k\le m$, then
$ h(x)\le E(x)$ for all $x\in \Delta_n$.
\end{thm}

\begin{proof}
Let $x\in \Delta_n$.  Also let $(X,\mu,\bpi)$ be an $S$-ranked
probability measure and $\la I_i\ra_{i=1}^\infty$ a disjoint sequence
in $\bpi$ such that
$$
\sum_{i=1}^\infty \mu(I_i)=1\quad \text{and}\quad \sum_{i\in
N_k}\mu(I_i)=x_k
$$
where $\N$ is partitioned by $N_0,N_1\cd N_n$ and
$x=\sum_{k=0}^nx_ke_k$. If $A\in \sigma\{I_i: i=1,2,\dots\}$ (the
$\sigma$-algebra generated by $\{I_i:i=1,2,\dots\}$), i.e.\ 
$A=\cup\{I_i: I_i\subseteq A\}$, we define (for $A\ne\nothing$, so
that $\mu(A)>0$) 
$$
x_A:=\frac{1}{\mu(A)}\sum_{k=0}^n
\bigg(\sum_{i\in N_k,\ I_i\subseteq A}\mu(I_i)\bigg)e_k. 
$$
Then the map $A\mapsto \mu(A)x_A$ is a vector measure on $\sigma\{I_i:
i=1,2,\dots\}$.  Note that $x_X=x$ (as $X=\bigcup_{i=0}^\infty I_i$
except for a set of $\mu$-measure zero so that $X\in \sigma\{I_i:
i=1,2,\dots\}$).  For each $m=1,2,3,\dots$ let
$$
A_m:=\bigcup_{r_\mu(I_i)\le m}I_i,\quad \text{and}\quad 
\mathcal{R}_m:=\{ J\in \pi_m: J\cap A_m=\nothing\}.
$$
Note that if $r_\mu(I_i)>m$, then $I_i\subseteq J$ for some $J\in
\mathcal{R}_m$.  Since 
$$
\sum_{r_\mu(I_i)>m}\mu(I_i)=1-\mu(A_m)=\sum_{J\in \mathcal{R}_m}\mu(J)
$$
each $J\in \mathcal{R}_M$ is (except for a set of $\mu$-measure zero)
a disjoint union of countable many sets $I_i$ with $r_\mu(I_i)>m$ so
that $J\in \sigma\{I_i: i=1,2,\dots\}$.  We require the following
lemma to complete the proof.

\begin{lemma}\label{lemma1}
With $h$ as in the statement of Theorem~\ref{E-extr} 
\begin{equation}\label{tail-ineq}
 h(x)\le 
\sum_{r_\mu(I_i)\le m}r_\mu(I_i)\mu(I_i)
+\sum_{J\in\mathcal{R}_m}[m+ h(x_J)]\mu(J).
\end{equation}
\end{lemma}

Before proving the lemma we show that it implies the theorem.  As
$ h$ is bounded there is an $M$ so that $ h(x)\le M$ for all $x\in
\Delta_n$.  Therefore by the lemma

\begin{align*}
 h(x)&\le \sum_{r_\mu(I_i)\le m}r_\mu(I_i)\mu(I_i)
	+\sum_{J\in \mathcal{R}_m} m\mu(J) 
	+M\sum_{J\in \mathcal{R}_m}\mu(J)\\
&= \sum_{r_\mu(I_i)\le m}r_\mu(I_i)\mu(I_i)
	+m\sum_{r_\mu(I_i)>m} \mu(I_i) 
	+M\sum_{r_\mu(I_i)>m}\mu(I_i)\\
&\le \sum_{r_\mu(I_i)\le m}r_\mu(I_i)\mu(I_i)
	+\sum_{r_\mu(I_i)>m}r_\mu(I_i) \mu(I_i) 
	+M\sum_{r_\mu(I_i)>m}\mu(I_i)\\
&= \sum_{i=1}^\infty r_\mu(I_i)\mu(I_i)
	+M\sum_{r_\mu(I_i)>m}\mu(I_i).
\end{align*}
Since $\lim_{m\to\infty}\sum_{r_\mu(I_i)>m}\mu(I_i)=0$ this yields
$ h(x)\le \sum_{i=1}^\infty\mu(I_i)r_\mu(I_i)$.  Taking the infimum
over $\mu$ gives $ h(x)\le E(x)$ and completes the proof of
Theorem~\ref{E-extr}. 
\end{proof}

\begin{proof}[Proof of Lemma~\ref{lemma1}]
The proof is by induction on $m$.  The base case is $m=0$ which
amounts to $ h(x)\le [0+ h(x_X)]\mu(X)$, which is in fact an
equality. Now assume for some $m\ge 0$ that the
inequality\eq{tail-ineq} holds.  Consider $J\in \mathcal{R}_m$.  Then
$J$ divides into sets $J_0,J_1\cd J_N\in \pi_{m+1}$ such that
$$
\(\frac{\mu(J_0)}{\mu(J)}\cd \frac{\mu(J_N)}{\mu(J)}\)\in S.
$$
Since $\sum_{i=0}^N\mu(J_i)x_{J_i}=\mu(J)x_J$ the $S$-almost convexity
of $ h$ implies
$$
 h(x_J)\le 1 + \sum_{i=0}^N\frac{\mu(J_i)}{\mu(J)} h(x_{J_i}).
$$
Multiplying this by $\mu(J)$ 
\begin{align*}
\mu(J) h(x_J)&\le \mu(J)+\sum_{i=0}^N\mu(J_i) h(x_{J_i})\\
&=\sum_{i=0}^N\mu(J_i)+\sum_{i=0}^N\mu(J_i) h(x_{J_i})\\
&=\sum_{i=0}^N[1+ h(x_{J_i})]\mu(J_i).
\end{align*}
If we let $\mathcal{S}_m=\{J\in \pi_{m+1}: J\cap A_m=\nothing\}$
and apply the
above to each $J\in \mathcal{R}_m$
$$
\sum_{J\in \mathcal{R}_m}[m+ h(x_J)]\mu(J)\le \sum_{J\in
\mathcal{S}_m}[m+1+ h(x_J)]\mu(J).
$$
If $J\in \mathcal{S}_m$ and $J=I_i$ for some $i$ then
$x_J=x_{I_i}=e_k$, where $i\in N_k$.  Thus the term for $J$ satisfies
$$
[m+1+ h(x_J)]\mu(J)=[m+1+ h(e_k)]\mu(I_i)\le(m+1)\mu(I_i)=r_\mu(I_i)\mu(I_i)
$$
since $I_i\in \pi_{m+1}$ and $h(e_k)\le 0$.
Now $\{J\in \mathcal{S}_m: J\ne I_i \text{ for any }
i\}=\mathcal{R}_{m+1}$.  Thus
\begin{align*}
 h(x)&\le \sum_{r_\mu(I_i)\le m}r_\mu(I_i)\mu(I_i)
	+\sum_{J\in\mathcal{R}_m}[m+ h(x_J)]\mu(J)\\
&\le \sum_{r_\mu(I_i)\le m}r_\mu(I_i)\mu(I_i)
	+\sum_{J\in\mathcal{S}_m}[m+1+ h(x_J)]\mu(J)\\
&\le\sum_{r_\mu(I_i)\le m}r_\mu(I_i)\mu(I_i)
	+\sum_{r_\mu(I_i)=m+1}r_\mu(I_i)\mu(I_i)\\
&\quad\qquad
	+\sum_{J\in\mathcal{R}_{m+1}}[m+1+ h(x_J)]\mu(J)\\
&=\sum_{r_\mu(I_i)\le m+1}r_\mu(I_i)\mu(I_i)
	+\sum_{J\in\mathcal{R}_{m+1}}[m+1+ h(x_J)]\mu(J).
\end{align*}
This closes the induction and completes the proof of the lemma.
\end{proof}

\subsection{Bounds for $S$-almost convex functions and the sharp
constants in stability theorems of Hyers-Ulam
type.}\label{sec:bds-Ulam} Let $S\subseteq
\bigcup_{m=1}^\infty\Delta_m$ and assume  that $S$ contains at
least one point that is not a vertex, that is a point $(t_0\cd t_m)$
with $\max_it_i<1$. Then, letting $E_S^{\Delta_n}\cn \Delta_n\to \R$
be as in Definition~\ref{def:E}, set
\begin{equation}\label{def:kappa}
\kappa_S(n):=\sup_{x\in \Delta_n}E_S^{\Delta_n}(x).
\end{equation}
By Proposition~\ref{prop1} the number $\kappa_S(n)$ is finite and we
will show that it is given by~\eqref{kappa(m)}.  The function
$E_S^{\Delta_n}$ and the number $\kappa_S(n)$ are extremal in several
analytic and geometric inequalities involving $S$-almost convex
functions and sets.  An example of this is the sharp form of the
Hyers-Ulam stability theorem (Theorem~\ref{stability}) in which
$\kappa_S(n)$ is the best constant and the example showing that this
is the case is the function~$E^{\Delta_n}_S$.  The exact value of
$\kappa_S(n)$ for some natural choices of $S$ are given in later
sections.  As a preliminary to Theorem~\ref{stability} we show that
$S$-almost convex functions with minimal regularity (Borel
measurability) are locally bounded so that Theorem~\ref{E-extr} can be
applied.

Recall that in a metric space the Borel sets are the members of the
$\sigma$-algebra generated by the open sets and if $f\cn X\to Y$ is a
function between metric spaces then it is Borel measurable iff
$f^{-1}[U]$ is a Borel subset of $X$ for every open subset $U$ of $Y$.

\begin{prop}\label{Borel-bded}
Assume that $S$ has at least one point that is not a vertex.  Let $
h\cn \Delta_n\to \R$ be a Borel measurable  $S$-almost convex function.
Then
\begin{align*}
 h(x)&\le E_S^{\Delta_n}(x)+x_0 h(e_0)+\cdots +x_n h(e_n)\\
&\le \kappa_S(n)+x_0 h(e_0)+\cdots +x_n h(e_n).
\end{align*}
\end{prop}

\begin{proof}
By replacing $ h$ by $x\mapsto  h(x)-(x_0 h(e_0)+\cdots
+x_n h(e_n))$, which will still be $S$-almost convex, we may assume
that $ h(e_i)=0$ for $0\le i\le n$.  If $ h$ is bounded then
$ h\le E_S^{\Delta_n}$ by Theorem~\ref{E-extr}.  So to finish the
proof it is enough to show that $ h$ is bounded.  In doing this we
can use Proposition~\ref{S-part} and note that there are
$\alpha,\beta>0$ with $\alpha+\beta=1$ so that if
$S_2=\{\alpha,\beta\}$ then $ h$ is $S_2$-almost convex.  (To be a
bit more precise let $(t_0\cd t_m)\in S$ with $\max t_i<1$ and then
the choice $\alpha=\max_it_i$ and $\beta=1-\alpha$ works.)  

With this choice of $S_2$ we now prove by induction on $n$ that if $
h\cn \Delta_n\to \R$ is $S_2$-almost convex and vanishes on the
vertices of $\Delta_n$ then $ h\le \kappa_{S_2}(n)$.  The base case is
$n=1$.  Then as a Borel measurable function is Lebesgue measurable
Corollary~\ref{1D-bded} implies $ h$ is bounded.  But then
Theorem~\ref{E-extr} implies $ h(x)\le E_{S_s}^{\Delta_1}(x)\le
\kappa_{S_2}(1)$.

For the induction step let $ h\cn \Delta_{n}\to \R$ be $S_2$-almost
convex and suppose $h$ vanishes on the vertices of $\Delta_n$.  Let
$g\cn \Delta_{n-1}\to \R$ be the function $g(y_0\cd y_{n-1})= h(y_0\cd
y_{n-1},0)$.  Then $g$ is $S_2$-almost convex, vanishes on the
vertices of $\Delta_{n-1}$ and is Borel measurable.  Therefore by the
induction hypothesis $g\le \kappa_{S_2}(n-1)$.  Let $y\in
\Delta_{n-1}$ and consider the function $\widetilde{ h}\cn [0,1]\to
\R$ given by
$$
\widetilde{ h}(t)= h((1-t)(y,0)+te_n)-(1-t) h(y,0).
$$
Then this is $S_2$-almost convex on $[0,1]$ and is Borel measurable.
Therefore another application of Corollary~\ref{1D-bded} implies that
$\widetilde{ h}$ bounded and as $\widetilde{ h}$ vanishes at the
endpoints of $[0,1]$ we have that $\widetilde{ h}(t)\le
\kappa_{S_2}(1)$.  This implies
\begin{align*}
 h((1-t)(y,0)+te_n)&=\widetilde{ h}(t)+(1-t) h(y,0)
=\widetilde{ h}(t)+(1-t)g(y)\\
&\le \kappa_{S_2}(1)+(1-t)\kappa_{S_2}(n-1)\\
&\le \kappa_{S_2}(1)+\kappa_{S_2}(n-1).
\end{align*}
But every $x\in \Delta_n$ can be expressed as $x=(1-t)(y,0)+te_n$ for
some $y\in \Delta_{n-1}$ and some $t\in [0,1]$.  Therefore $ h$ is
bounded on $\Delta_n$.  Then Theorem~\ref{E-extr} implies $ h(x)\le
E_{S_2}(x)\le \kappa_{S_2}(n)$.  This closes the induction and
completes the proof.
\end{proof}

\begin{thm}\label{gen-bds}
Let $U$ be a convex set in a normed vector space and let $h\cn U\to
\R$ be an $S$-almost convex function which is bounded above on compact
subsets of $U$.   Assume that $S$ contains at least one point which is
not a vertex.  Then for any $x_0\cd x_n\in U$ the inequalities
\begin{align}
h(t_0x_0+\cdots+t_nx_n)&\le
E_{S}^{\Delta_n}(t)+t_0h(x_0)+\cdots+t_nh(x_n)\nn\\
&\le \kappa_S(n)+t_0h(x_0)+\cdots+t_nh(x_n)\label{gen-bd}
\end{align}
hold for all $t=(t_0\cd t_n)\in \Delta_n$.  If $U$ is compact,
$n$-dimensional and $V$ is the set of extreme points of $U$ then
\begin{equation}\label{V-U-bd}
\sup_{x\in U}h(x)\le \kappa_S(n)+\sup_{v\in V}h(v).
\end{equation}
\end{thm}

\begin{proof}
Let $f\cn \Delta_n\to \R$ be given by $f(t)=h(t_0x_0+\cdots+t_nx_n)
-(t_0h(x_0)+\cdots+t_nh(x_n))$.  Then $f$ is $S$-almost convex,
bounded (as $h$ is bounded on the convex hull of $\{x_0\cd x_n\}$ as
it is compact) and vanishes on the vertices of $\Delta_n$.  Therefore
by Theorem~\ref{E-extr} $f(t)\le E_{S}^{\Delta_n}(t)\le\kappa_S(n)$
which implies\eq{gen-bd}.

If $U$ is compact and $n$ dimensional with extreme points $V$, then $U$
is the convex hull of $V$.  By Carath\'eodory's Theorem for any $x\in
U$ there are $x_0\cd x_n\in V$ and $t=(t_0\cd t_n)$ so that
$x=t_0x_0+\cdots+t_nx_n$ which, along with\eq{gen-bd}, implies\eq{V-U-bd}.
\end{proof}

We can now give the sharp version of the Hyers-Ulam stability theorem
for $S$-almost convex functions.

\begin{thm}\label{stability}
Let $S\subseteq \bigcup_{m=1}^\infty\Delta_m$ so that $S$ contains at
least one point that is not a vertex.  Assume that $U\subseteq \R^n$,
$\e>0$, and that $ h\cn U\to \R$ is bounded above on compact subsets
of $U$ and satisfies
\begin{equation}\label{e-ineq}
 h(t_0x_0+\cdots+ t_mx_m)\le \e +t_0 h(x_0)+\cdots +t_m h(x_0)
\end{equation}
for all $t=(t_0\cd t_m)\in S$ and points $x_0\cd x_m\in U$.  Then
there exist convex functions $g,g_0\cn U\to \R$ such that
$$
 h(x)\le g(x)\le  h(x)+\kappa_S(n)\e\qquad \text{and}\qquad
| h-g_0(x)|\le \frac{\kappa_S(n)}{2}\e
$$
for all $x\in U$.  The constant $\kappa_S(n)$ is the best constant in
these inequalities.
\end{thm}

\begin{remark}
Note that if $ h$ satisfies\eq{e-ineq} then $\e^{-1}h$ is $S$-almost
convex.  Therefore, by Theorem~\ref{leb-bded}, if $U$ is open and
$ h$ is Lebesgue measurable then $ h$ will automatically be
bounded on compact subsets of $U$.  Likewise if $U$ is a Borel set and
$h$ is Borel measurable then by Proposition~\ref{Borel-bded} $h$ will
be bounded above on the convex hull of any finite number of points and
this is enough for the proof of the theorem.\qed
\end{remark}

\begin{proof} In the special case that
$S=\{(1/2,1/2)\}\subset\Delta_1$ a proof, based on ideas of Hyers and
Ulam~\cite[p.~823]{Hyers-Ulam:convex} and
Cholewa~\cite[pp.~81--82]{Cholewa:stability}, can be found
in~\cite[pp.~29-30]{Dilworth-Howard-Roberts:1}.  As the details in the
present case are identical we omit the proof.
\end{proof}


\section{General results when $S$ is compact.}
\label{sec:S-cmpt}

We now assume that $S\subseteq  \bigcup_{m=1}^\infty\Delta_m$ is
compact.  By Remark~\ref{S-compact} this implies that $S$ is of finite
type.  Therefore by Proposition~\ref{1Delta} there is no loss in
generality in assuming that $S\subseteq\Delta_m$ for some $m$.  

\subsection{Mean value and semi-continuity properties.}\label{sec:mean-lsc}
Let $K\subset\R^n$ be a compact convex set and let $V$ be the set of
extreme points of $K$. If $\phi\cn V\to \R$ is a function, then $h\cn
K\to \R$ has {\bi extreme values equal to $\phi$\/} iff
$h\big|_V=\phi$.  Two functions $g,f\cn K\to \R$ {\bi have the same
extreme values\/} iff they agree on $V$.  If $\phi\cn V\to \R$ is a
bounded function and  $S\subseteq \Delta_m$ then let
$\mathcal{B}_S(K,\phi)$ be the set of bounded $S$-almost convex
functions $h\cn K\to \R$ so that $h\big|_V\le \phi$.  Then the {\bi
extremal $S$-almost convex function with extreme values $\phi$\/} is 
$$
E_{S,K,\phi}(x):=\sup_{h\in \mathcal{B}_S(K,\phi)}h(x).
$$
If $S$ contains at least one point which is not a vertex, then
Theorem~\ref{gen-bds} implies that $E_{S,K,\phi}$ is finite valued and
in fact $E_{S,K,\phi}(x)\le \sup_{v\in V}\phi(v)+\kappa_S(n)$.  As the
pointwise supremum of $S$-almost convex functions is $S$-almost convex,
the function $E_{S,K,\phi}$ is the pointwise largest $S$-almost convex
function with $E_{S,K,\phi}(v)\le \phi(v)$ on $V$.

If $K\subset\R^n$ is a compact convex set and $V$ is the set of
extreme points of $K$ then for any function $h\cn K\to \R$ define
$\mean_Sh\cn K\to \R^n$ by
$$
\mean_Sh(x)=\begin{cases}h(x),& x\in V;\\
	\inf\left\{\displaystyle1+\sum_{i=0}^mt_ih(y_i)\ :\  t\in S,\ 
	x=\sum_{i=0}^mt_iy_i\right\},& x\in K\setminus V\end{cases}
$$
where it is assumed that $y_0\cd y_m\in K$.  We can then define
$S$-almost convex functions in terms of this operator by the
following, for any bounded function $f\cn K\to \R$,
$$
f\le \mean_Sf\quad \iff \quad \text{$f$ is $S$-almost convex.}
$$
This operator satisfies a maximum principle and can be used to prove
that extremal $S$-almost convex functions are lower semi-continuous.

\begin{thm}\label{max-prin}
Let $K\subset\R^n$ be a compact convex set with extreme points
$V$. Assume that $S\subset \Delta_m$ is compact and has at least one
point which is not a vertex.  Let $f,F\cn K\to \R$ be bounded
functions so that $\mean_Sf\le f$ and $F$ is $S$-almost convex (that
is $\mean_SF\ge F$) then
\begin{equation}\label{f-max}
\sup_{x\in K}(F(x)-f(x))=\sup_{v\in V}(F(v)-f(v))
\end{equation}
and if $L$ is the lower semi-continuous envelope of
$f$,
\begin{equation}\label{L-max}
L(x):= \min\{f(x),\liminf_{y\to x} f(y)\}
\end{equation}
then 
$$
\sup_{x\in K}\(F(x)-L(x)\)=\sup_{v\in V}\(F(v)-L(v)\).
$$
\end{thm}

\begin{remark}
The proof here follows the basic outline of the proof of corresponding
result, \cite[Theorem~2.8 p.9]{Dilworth-Howard-Roberts:1}, in the case
$S=\{(1/2/1/2)\}\subset\Delta_1$.  However the technical details are
trickier in the case when $S$ is infinite.  But most of the rest of
the results of~\cite[Section~2.2]{Dilworth-Howard-Roberts:1} go
through with only minor changes to the proofs.\qed
\end{remark}

\begin{proof}
The proofs of\eq{f-max} and\eq{L-max} are similar, with the proof
of\eq{f-max} being the simpler of the two, so we will give the details
in the proof of\eq{L-max}.  The inequality $f\ge \mean_Sf$ implies for
$x\notin V$ and any $\e>0$ there is a $t=(t_0\cd t_m)\in S$ and
$y_0\cd y_m\in K$ such that
\begin{equation}\label{f-ineq}
x=\sum_{i=0}^mt_i y_i,\qquad f(x)\ge 1-\e+\sum_{i=0}^mt_if(y_i).
\end{equation}
As $f$ and $F$ are bounded we can assume, by adding appropriate
positive constants to $f$ and $F$, that $1\le f\le F\le M$ for
some $M>1$.  Set
$$
\omega(x):=F(x)-L(x),\quad \delta:=\sup_{x\in K}\omega(x).
$$
We need to show that $\sup_{v\in V}\omega(v)\ge \delta$ (as
$\sup_{v\in V}\omega(v)\le \delta$ is clear).  We may assume that
$\delta>0$, for if $\delta=0$ then $F=L$ and there is nothing to
prove. 

\begin{lemma}\label{step-down}
Let $w_0\in K$, but $w_0\notin V$ and assume for some $\e>0$  that 
$$
(1-\e)\delta \le \omega(w_0).
$$
Then there is a $w_1\in K$ so that
$$
\(1-(m+1)(2M-1)\e\)\delta \le \omega(w_1)\quad \text{and} \quad
L(w_1)\le L(w_0)-\frac12.
$$
\end{lemma}

We now prove Theorem~\ref{max-prin} from the lemma.  Let $\e>0$. We
now choose a finite sequence $w_0, w_1\cd w_k$ with $k\le 2M$ as
follows.  From the definition of $\delta$ there is a $w_0\in K$ with
$(1-\e)\delta \le \omega(w_0)$.  If $w_0\in V$ we stop.  If $w_0\notin
V$, then by the lemma, there is a $w_1\in K$ with
$(1-(m+1)(2M-1)\e)\delta \le \omega(w_1)$ and $L(w_1)\le L(w_0)-1/2$.
If $w_1\in V$ then stop, otherwise use the lemma (with $w_1$ replacing
$w_0$ and $(m+1)(2M-1)\e$ replacing $\e$) to get a $w_2$ with
$(1-((m+1)(2M-1))^2\e)\delta \le \delta$.  If $w_2\in V$, stop.  If
$w_2\notin V$ then we continue to use the lemma to get $w_0,w_1\cd
w_k$ with
$$
\(1-\((m+1)(2M-1)\)^j\e\)\delta\le \omega(w_j)\quad \text{and} \quad
L(w_j)\le L(w_{j-1})-\frac12
$$
for $1\le j\le k$.  This implies that $L(w_k)\le L(w_0)-k/2\le M-2/k$.
But as $L\ge1$ this process must terminate for some $k\le 2M$ with
$w_k\in V$.  Then
\begin{align*}
\sup_{v\in V}\omega(v)&\ge \omega(w_k)\ge
\(1-\((m+1)(2M-1)\)^{k}\e\)\delta\\
&\ge \(1-\((m+1)(2M-1)\)^{2M}\e\)\delta.
\end{align*}
Letting $\e\searrow0$ in this implies $\sup_{v\in V}\omega(v)\ge
\delta$ which completes the proof.
\end{proof}

\begin{proof}[Proof of Lemma~\ref{step-down}]
Let $w_0$ be as in the statement of the lemma.  From the definition of
$L$ there is a sequence $\la x(s)\ra_{s=1}^\infty\subset K$ so that
$x(s)\to w_0$ and 
$f(x(s))\to L(w_0)$.  By\eq{f-ineq} there is a sequence 
$\la t(s)\ra_{s=1}^\infty=\la (t_0(s)\cd t_m(s))\ra_{s=1}^\infty\subseteq S$
and sequences $\la y_0(s)\ra_{s=0}^\infty\cd \la
y_m(s)\ra_{s=0}^\infty \subseteq K$ so that (replacing $\la
x(s)\ra_{s=1}^\infty$ by the appropriate subsequence).
$$
f(x(s))-\(1+\sum_{i=0}^m t_i(s)f(y_i(s))\)\xrightarrow{s\to\infty} C\ge 0
$$
for some non-negative real number $C$.  By compactness of $S$ and $K$
we can assume, by possibly going to a subsequence, that $t(s)\to t\in
S$ and $y_i(s)\to y_i\in K$ and that $f(y_i(s))\to A_i$ for some $t\in
S$, $y_0\cd y_m\in S$ and $A_i\in \R$.  Then $w_0=\sum_{i=0}^mt_iy_i$
and from the definition of $L$, $L(y_i)\le \lim_{s\to
\infty}f(y_i(s))=A_i$.  Therefore
\begin{equation}\label{Lw0} 
\lim_{s\to\infty}f(x(s))=
L(w_0)=C+1+\sum_{i=0}^m t_iA_i\ge 1+\sum_{i=0}^mt_iL(y_i).
\end{equation}
This is turn implies that
\begin{equation}\label{Fw0big}
F(w_0)=\omega(w_0)+L(w_0)\ge \omega(w_0)+1+ \sum_{i=0}^mt_iL(y_i).
\end{equation}
Because $F$ is $S$-almost convex,
\begin{equation}\label{Fw0small}
F(w_0)\le 1+\sum_{i=0}^mt_iF(y_i)=1+\sum_{i=0}^mt_iL(y_i)
	+\sum_{i=0}^mt_i\omega(y_i).
\end{equation}
Combining\eq{Fw0big} and\eq{Fw0small} yields
\begin{equation}\label{omega-avg}
\omega(w_0)\le \sum_{i=0}^mt_i\omega(y_i).
\end{equation}

We now claim there is an $i_0$ so that 
\begin{equation}\label{L-t-claim}
t_{i_0}\ge \frac1{(m+1)(2M-1)},\quad L(y_{i_0})\le L(w_0)-\frac12.
\end{equation}
To see this partition $\{0,1\cd m\}$ into two sets $I_1$ and $I_2$
where $I_1:=\{i: t_i< 1/((m+1)(2M-1))\}$ and $I_2:=\{i: t_i\ge
1/((m+1)(2M-1))\}=\{0\cd m\} \setminus I_1$.  Note that as $M>1$ we
have 
$$
\sum_{i\in I_1}t_i\le (m+1)/((m+1)(2M-1))=1/(2M-1)<1/2
$$ 
so that $I_2\ne \nothing$.  For $i\in I_2$ let $\alpha_i=(\sum_{i\in
I_2}t_i)^{-1}t_i$.  Then $\sum_{i\in I_2}\alpha_i=1$.  Using\eq{Lw0},
\begin{align*}
\sum_{i\in I_2}\alpha_iL(y_i)
&=\bigg(\sum_{i\in I_2}t_i\bigg)^{-1}\sum_{i\in I_2}t_iL(y_i)
\le\bigg(\sum_{i\in I_2}t_i\bigg)^{-1}\sum_{i=0}^mt_iL(y_i)\\
&\le \bigg(\sum_{i\in I_2}t_i\bigg)^{-1}(L(w_0)-1).
\end{align*}
We have already seen that $1-\sum_{i\in I_2}t_i=\sum_{i\in I_1}t_i\le
1/(2M-1)$ and therefore $\sum_{i\in I_2}t_i\ge 1-1/(2M-1)=(M-1)/(M-1/2)$.
Thus
\begin{align*}
\sum_{i\in I_2}\alpha_iL(y_i)&\le \frac{M-1/2}{M-1} (L(w_0)-1)\\
&\le \frac{L(w_0)-1/2}{L(w_0)-1}(L(w_0)-1)\\
&=L(w_0)-\frac12
\end{align*}
where we have used that $L(w_0)\le M$ and that $(M-1/2)/(M-1)$ is
decreasing for $M>1$.  As $\sum_{i\in I_2}\alpha_i=1$ this implies
there is at least one $i_0\in I_2$ with $L(y_{i_0})\le L(w_0)-1/2$.
For this $i_0$ the claim\eq{L-t-claim} holds.

Letting $i_0$ be so that\eq{L-t-claim} holds and using that 
$(1-\e)\delta\le \omega(w_0)$, and that
$\omega(y_i)\le \delta$ for all $i$ in\eq{omega-avg}, we have
$$
\(1-\e \)\delta
\le \omega(w_0) \le  \sum_{i=0}^mt_i\omega(y_i) 
\le t_{i_0}\omega(y_{i_0})+(1-t_{i_0})\delta.
$$
This implies
$$
\(1- t_{i_0}^{-1}\e\)\delta\le \omega(y_{i_0}).
$$
As $t_{i_0}\ge 1/((m+1)(2M-1))$ this gives
$$
\(1-\((m+1)(2M-1)\)\e\)\delta\le \omega(y_{i_0}).
$$
Letting $w_1=y_{i_0}$  completes the proof of the lemma.
\end{proof}

\begin{thm}\label{extr-lsc}
Let $K\subset\R^n$ be a compact convex set with extreme points $V$.
Assume that $\phi\cn V\to \R$ is uniformly continuous.  Let
$S\subseteq \Delta_m$ be compact and contain at least one point that
is not a vertex.  Then the extremal $S$-almost convex function
$E_{S,K,\phi}$ is lower semi-continuous and satisfies
$E_{S,K,\phi}\big|_V=\phi$.
\end{thm}

\begin{proof}
This can be derived from Theorem~\ref{max-prin} in the same way that
\cite[Theorem~2.12 p.~13]{Dilworth-Howard-Roberts:1} is derived from 
\cite[Theorem~2.8 p.~9]{Dilworth-Howard-Roberts:1}.
\end{proof}

\subsection{Simplifications in the  construction of
$E^{\Delta_n}_S$ when $S$ is compact.}
\label{sec:simple}
One complication in Definition~\ref{def:E} is that the infimum is
taken over a collection of measures that are not all defined on the
same measure space.  When $S\subseteq \Delta_m$ it is possible to have
all the measures involved defined on the same space.

Suppose $S\subseteq \Delta_m$.  We may regard each $S$-ranked
probability measure as a (Borel) probability measure on $X=[m]^\N$,
with $[m]=\{0,1\cd m\}$.
Let $\prob(X)$ be the space of probability measures on $X$.  Then
$\prob(X)\subset C(X)^*$ and in the weak$^*$ topology $\prob(X)$ is
compact and metrizable (as $C(X)$ is separable).  We let
$$
\prob_S(X):=\{ \mu\in \prob(X): \text{$\mu$ is $S$-ranked}\}.
$$
Then every $\mu\in \prob_S(X)$ has $\pi_{j}(\mu)=\pi_j$ given by
$$
I\in \pi_j \quad \iff\quad \left\{\begin{array}{l}
 \text{for some } (i_1\cd i_j)\in [m]^j, \\
	I=\{x\in X: x(1)=i_1\cd x(j)=i_j\}.\end{array}\right.
$$
or what is the same thing $I\in \pi_j$ if and only if
$I=\{i_1\}\times \{i_2\}\times \{i_j\} \times X_j$ where
$X_j=\prod_{i=j+1}^\infty Y_i$ with $Y_i=[m]$ for all $i$.  Since each
$\mu \in \prob_S(X)$ has the same sequence $\bpi=\la \pi_j\ra$, we
let $r(I)=r_\mu(I)$ which is defined independently of the choice of
$\mu\in \prob_S(X)$.  Let $\pi =\bigcup_{j=1}^\infty \pi_j$.  

Finally note that if $\alg_j=\alg(\pi_j)$ and $A\in \alg_j$, then $A$
is a clopen (i.e.\ both open and closed) set in $X$.  Consequently
$\boldsymbol{1}_A\in C(X)$.  In this case we have $\alg_j=\alg(\pi_j)$
and thus the function $\mu\mapsto \mu(A)=\int\boldsymbol{1}_A\,d\mu$
is continuous on $\prob(X)$ and thus on $\prob_S(X)$.\qed

\begin{prop}\label{prop3}
With this notation, if $S\subset \Delta_m$ is
closed, then $\prob_S(X)$ is closed in $\prob(X)$ and thus is weak$^*$
compact.
\end{prop}

\begin{proof}
Notice that if $\mu\in \prob(X)$, then $\mu\in \prob_S(X)$ if and only
if for every $I\in \pi$, there exists $(t_0,t_1\cd t_m)\in S$ such
that
$$
\mu(I)=\sum_{i=0}^mt_i\mu(I_i)
$$
where $I\in \pi_j$ and $I$ is the disjoint union of $I_0,I_1\cd I_m\in
\pi_{j+1}$.  Let ${t}=(t_0,t_1\cd t_m)\in S$ and define a function
$ h_{I,{t}}\cn \prob(X)\to \R$ by
$$
 h_{I,{t}}(\mu)=\mu(I)-\sum_{i=0}^mt_i\mu(I_i).
$$
Then this is continuous on $\prob(X)$.  Let
$$
\Lambda_I:=\bigcup_{t\in S}  h^{-1}_{I,t}[\{0\}]=
\{\mu\in \prob(X): \mu(I)=\sum_{i=0}^mt_i\mu(I_i)\text{ for some }
t\in S\}.
$$
Then $\prob_S=\bigcap_{I\in \pi}\Lambda_I$.  As an intersection of
closed sets is closed, to finish the proof it is enough to show that
each $\Lambda_I$ is closed.  Let $\mu_s\in \Lambda_I$ and suppose
$\mu_s\xrightarrow{\text{\it weak$^*$}} \mu$ in $\prob(X)$.  For each
$s=1,2,3,\dots$ there is a $t(s)=(t_0(s)\cd t_m(s))\in S$ such that
$\mu_s(I)=\sum_{i=0}^m t_i(s)\mu_s(I_i)$.  Since $S$ is compact, by
passing to a subsequence, if necessary, we may assume that $t(s)\to
t=(t_0\cd t_m)\in S$.  Thus
$$
\mu(I)=\lim_{s\to \infty}\mu_s(I)=\lim_{s\to
\infty}\sum_{i=0}^mt_i(s)\mu_s(I_i)=\sum_{i=0}^mt_i\mu(I_i).
$$
Therefore $\Lambda_I$ is closed.
\end{proof}

\begin{prop}\label{prop4}
Suppose that $S\subset\Delta_m$ is closed and that $S$ contains a
point that is not a vertex (so that by Proposition~\ref{prop1} $E=E_S$
is bounded).  Then
\begin{enumerate}
\item $E$ is lower semi-continuous,
\item If $x\in \Delta_n$, then there exists a $\mu\in \prob_S(X)$ and
a pairwise disjoint sequence $\la I_i\ra\in \pi$ such that
$$
\sum_{i=1}^\infty\mu(I_i)=1, \qquad x\mid \la \mu(I_i)\ra
$$
and
\begin{equation}\label{inf=min}
E(x)=\sum_{i=1}^\infty \mu(I_i)r(I_i).
\end{equation}
Thus the infimum that defines $E(x)$ is a minimum.
\end{enumerate}
\end{prop}

\begin{remark}\label{alt-lsc}
The lower semi-continuity of $E$ also follows from
Theorem~\ref{extr-lsc}, but we include another proof here both because
it is short and also to have a proof that is independent
of~\cite{Dilworth-Howard-Roberts:1}.\qed
\end{remark}

\begin{lemma}\label{lemma2}
Suppose that $S$ is a closed subset of $\Delta_m$.  Further suppose 
\begin{enumerate}
\item\label{1} $\la x(s)\ra_{s=1}^\infty$ is a sequence in 
$\Delta_n$ with $x(s)\to x\in \Delta_n$,
\item\label{2} $\la \mu_s\ra_{s=1}^\infty$ is a sequence in $\prob_S(X)$ with
$\mu_s\xrightarrow{\text{\it weak$^*$}}\mu \in \prob_S(X)$,
\item\label{3} For all $s\in \N$, there exists a disjoint sequence $\la
I_{j\,s}\ra_{j=1}^\infty\subset \{\nothing\}\cup\bigcup_{l=0}^\infty
\pi_l$ such that $\sum_{j=1}^\infty\mu_s(I_{j\,s}) =1$,
\item\label{4} $x(s)\mid \la \mu(I_{j\,s})\ra_{j=1}^\infty$, and
\item\label{5} There is an $M>0$ so that for all $s\in \N$
$$
M_s:=\sum_{j=1}^\infty r(I_{j\,s})\mu_s(I_{j\,s}) \le M.
$$
\end{enumerate}
Then there exists a disjoint sequence $\la
I_j\ra_{j=1}^\infty\subset\{\nothing\}\cup\bigcup_{l=0}^\infty \pi_l$ such
that
\begin{enumerate}
\item[i.]\label{i} $\displaystyle\sum_{j=1}^\infty \mu(I_j)=1$
\item[ii.]\label{ii} $x\mid \la \mu(I_j)\ra_{j=1}^\infty$
\item[iii.] \label{iii} $\displaystyle\sum_{j=1}^\infty r(I_j)\mu(I_j)\le \limsup_{s\to \infty} M_s$.
\end{enumerate}
\end{lemma}

\begin{proof}
First we select a subsequence $\la x_a\ra_{a\in F}$ of $\la
x_s\ra_{s=1}^\infty$ for some infinite $F\subseteq \N$ by first
choosing sets $F_j(k)\subseteq\N$ and $I_j(k)\in
\{\nothing\}\cup\bigcup_{l=0}^\infty \pi_l$ as follows: For each $s\in
\N$, we can use point~(\ref{4})  to partition the terms of
$\la I_{j\,s}\ra_{j=1}^\infty$ into $n+1$ sequences $\la
I_{j\,s}(0)\ra_{j=1}^\infty \cd \la I_{j\,s}(n)\ra_{j=1}^\infty$ where
$\sum_{j=1}^\infty
\mu(I_{j\,s}(k))=x_s(k)$ and $x_s=\sum_{k=0}^nx_s(k)e_k$.  We may
assume that for every $k\in\{0,1\cd n\}$ that $r(I_{1\,s}(k))\le
r(I_{2\,s}(k))\le \cdots$. If $\lim_{s\to
\infty}r(I_{1\,s}(0))=\infty$, let $F_1(0)=\N$ and $I_1(0)=\nothing$,
otherwise $\la r(I_{1\,s}(0))\ra_{s=1}^\infty$ is bounded for some
infinite set of $s\in\N$.  Since for any integer $L$, there are only
finitely many sets in $\pi$ of rank $\le L$, there is an $I_1(0)$ so
that $I_{1\,s}(0)=I_1(0)$ on an infinite subset $F_1(0)$ of $\N$.
Similarly choose $F_1(1)$ infinite in $F_1(0)$ and $I_1(1)$ such that
either $\lim_{s\to \infty}r(I_{1\,s}(1))=\infty$ and $I_1(1)=\nothing$
or $I_{1\,s}(1)=I_1(1)$ for all $s\in F_1(1)$.  Continue selecting
infinite sets $F_j(k)$ of $\N$ and $I_j(k)\in
\{\nothing\}\cup\bigcup_{l=0}^\infty \pi_l$ such that
$$
F_1(0)\supseteq F_1(1)\supseteq \cdots \supseteq F_1(n)\supseteq
F_2(0)\supseteq F_2(1)\supseteq	\cdots
$$
and either $\lim_{s\to\infty}r(I_{j\,s}(k))=\infty$ and
$I_{j}(k)=\nothing$ or $I_{j\,s}(k)=I_j(k)$ for
all  $s\in F_j(k)$.  The inequalities $r(I_{1\,s}(k))\le
r(I_{2\,s}(k))\le \cdots$ yield
$$
\lim_{s\to\infty}r(I_{j\,s}(k))=\infty\quad \text{implies}\quad
\lim_{s\to\infty}r(I_{j+1\,s}(k))=\infty,
$$
and therefore
$$
\lim_{s\to\infty}r(I_{j\,s}(k))=\infty\quad \text{implies}\quad
\nothing=I_{j}(k)=I_{j+1}(k)=I_{j+2}(k)=\cdots.
$$
Also the sets $\la I_j(k)\ra_{j=1}^\infty{}_{k=0}^{n}$ are pairwise
disjoint.

Now let $F$ be an infinite set in $\N$ such that each $F\setminus
F_j(k)$ is finite. Let $L\in \N$.  Assumption~(\ref{5}) implies
$$
(L+1)\sum_{r(I_{j\,s}(k))\ge L+1}\mu_s(I_{j\,s})
\le \sum_{r(I_{j\,s}(k))\ge L+1}r(I_{j\,s})\mu_s(I_{j\,s})\le M
$$
Thus for fixed $s$ and $k$ 
$$
\sum_{r(I_{j\,s}(k))\ge L+1}\mu_s(I_{j\,s})\le \frac{M}{L+1}
$$
and therefore
$$
\sum_{r(I_{j\,s}(k))\le L}\mu_s(I_{j\,s})\ge x_s(k)-\frac{M}{L+1}.
$$
Hence
\begin{align*}
\sum_{r(I_{j\,s}(k))\le L}\mu(I_j(k))&=
\lim_{\substack{s\in F\\s\to\infty}}\sum_{r(I_{j\,s}(k))\le
L}\mu_s(I_{j\,s})\\
&\ge\lim_{s\to\infty}x_s(k)-\frac{M}{L+1}\\
&=x(k)-\frac{M}{L+1}
\end{align*}
where $x=\sum_{k=0}^nx(k)e_k$.  It follows that
$$
\sum_{j=1}^\infty \mu(I_j(k))\ge x(k).
$$
But since the sets $ \la I_j(k)\ra$ are pairwise disjoint
$$
1\ge \sum_{k=0}^n\sum_{j=1}^\infty\mu(I_j(k))\ge \sum_{k=0}^nx(k)=1.
$$
But this implies that there must be equality for each $k\in \{0\cd
n\}$:
$$
\sum_{j=1}^\infty \mu(I_j(k))=x(k).
$$

Once again fix $L\in \N$.  For $s$ suitably large in $F$,
$I_{j\,s}(k)=I_j(k)$ if $r(I_j(k))\le L$.  Thus
\begin{align*}
\sum_{r(I_j(k))\le L} \mu(I_j(k))r(I_j(k))
&=\lim_{s\to\infty}\sum_{r(I_j(k))\le L} \mu(I_{j\,s}(k))r(I_{j\,s}(k))\\
&\le \limsup_{s\to\infty}M_s.
\end{align*}
(All the sums are finite so there is no problem in interchanging the
limit with the summation.)  Since this holds for all large $L\in \N$,
$$
\sum_{j,k} \mu(I_j(k))r(I_j(k))\le \limsup_{s\to\infty}M_s.
$$
Now splice the sequences $\la I_j(0)\ra_{j=1}^\infty,\la
I_j(1)\ra_{j=1}^\infty\cd \la I_j(n)\ra_{j=1}^\infty$ into a single
sequence $\la I_j\ra_{j=1}^\infty \subset
\{\nothing\}\cup\bigcup_{l=0}^\infty\pi_l$.  This sequence satisfies
the conclusion of the Lemma.
\end{proof}

\begin{proof}[Proof of Proposition~\ref{prop4}] We First show the
lower semi-continuity of $E$.  Suppose that $\la x(s)\ra_{s=1}^\infty$
is a sequence in $\Delta_n$ and that $x(s)\to x\in \Delta_n$.  Further
suppose that $\la E(x(s))\ra$ is convergent.  For each $s\in \N$,
select a measure $\mu_s\in \prob_S(X)$ and a sequence $\la
I_{j\,s}\ra_{j=1}^\infty$ in $\pi$ such that $\sum_{j=1}^\infty
\mu_s(I_{j\,s})=1$, $x(s)\mid \la \mu_s(I_{j\,s})\ra_{j=1}^\infty$, and
$M_s=\sum_{j=1}^\infty\mu_s(I_{j\,s})r(I_{j\,s})<E(x_s)+1/s$.  By
passing to a subsequence, if necessary, we may assume that
$\mu_s\xrightarrow{\text{\it weak$^*$}}\mu\in \prob_S(X)$.  By
Lemma~\ref{lemma2}, there is a sequence $\la I_i\ra_{i=1}^\infty$ in
$\pi$ so that $\sum_{i=1}^\infty \mu(I_i)r(I_i)=1$, $x\mid \la
\mu(I_i)\ra_{i=1}^\infty$ and
$$
E(x)\le \sum_{i=1}^\infty \mu(I_i)r(I_i)\le \limsup_{s\to\infty}
M_s=\lim_{s\to\infty}E(x_s). 
$$
Thus $E$ is lower semi-continuous.

We now show the second conclusion of Proposition~\ref{prop4}.  Let
$x\in \Delta_n$.  Select $\la \mu_s\ra_{s=1}^\infty $ a sequence in
$\prob_S(X)$ and for each $s$ choose a sequence $\la
I_{j\,s}\ra_{j=1}^\infty$ in $\pi$ such that
$\sum_{j=1}^\infty\mu_s(I_{j\,s})=1$, $x\mid \la
\mu_s(I_{j\,s})\ra_{j=1}^\infty$ and
$$
E(x)\le \sum_{j=1}^\infty \mu_s(I_{j\,s})r(I_{j\,s})<E(x)+\frac1s.
$$
By passing to a subsequence, if necessary,
$\mu_s\xrightarrow{\text{\it weak$^*$}}\mu$ for some
$\mu\in\prob_S(X)$.  Let $\la I_i\ra_{j=1}^\infty$ be the sequence
obtained by Lemma~\ref{lemma2}.  Then for the measure $\mu$ and the
sequence $\la I_j\ra_{j=1}^\infty$ the equality~\ref{inf=min} holds.
This completes the proof.
\end{proof}

\section{Explicit Calculation of $E_S^{\Delta_n}$ and $\kappa_S(n)$
when $S$ is the barycenter of $\Delta_m$.}

The most natural choices of $S$ are when $S$ is a entire simplex
$\Delta_m$ or $S$ is the barycenter of $\Delta_m$. We have treated the
case of $S=\Delta_m$ in a previous
paper~\cite{Dilworth-Howard-Roberts:3} by different methods.  Here we
compute $E^{\Delta_n}_S$ and $\kappa_S(n)$ in the case $S$ is the
barycenter of $\Delta_m$ based on the general theory above.  It will
simplify notation to let $B=m+1$.


We now assume that $S=\{(1/B\cd 1/B)\}\subset \Delta_{B-1}$. To give
$E_S^{\Delta_n}$ explicitly we need a little notation.  First for any
real number $x$ let $\fp{x}=x-\lfloor{x}\rfloor$ be the fractional
part of $x$ and define a function $H=H_B\cn \R\to \R$ from by
\begin{equation}\label{H-fp}
H_B(x)=\sum_{k=0}^\infty \frac{\fp{B^kx}}{B^k}.
\end{equation}
Note that this series is termwise dominated by the geometric series
$\sum_{k=0}^\infty1/B^k$ and therefore it is easy to deal with
computationally.

\begin{thm}\label{main-bary}
For $S=\{(1/B\cd 1/B)\}$ the function $E:=E_S\cn \Delta_n\to \R$ is given by
$$
E(x)=E(x_0,x_1\cd x_n)=H_B(x_0)+H_B(x_1)+\cdots+H_B(x_n)
$$
and the value of $\kappa_S(n)=\sup_{x\in\Delta_n}E(x)$ is 
$$
\kappa_S(n)=\lfloor\log_B n\rfloor+1+\frac{n}{(B-1)B^{\lfloor\log_Bn\rfloor}}.
$$
\end{thm}

Some values of $\kappa_S(n)$ for small values of $B$ and $n$ are given
in Table~\ref{kappa-table}.

\begin{table}[hb]
{\tiny
$$
\begin{array}{c|cccccccccc}
B\backslash n& 1&2&3&4&5&6&7&8&9&10\\ \hline
2&2.0000&3.0000&3.5000&4.0000&4.2500&4.5000&4.7500&5.0000&5.1250&5.2500 \\
3&1.5000&2.0000&2.5000&2.6667&2.8333&3.0000&3.1667&3.3333&3.5000&3.5556 \\
4&1.3333&1.6667&2.0000&2.3333&2.4167&2.5000&2.5833&2.6667&2.7500&2.8333 \\
5&1.2500&1.5000&1.7500&2.0000&2.2500&2.3000&2.3500&2.4000&2.4500&2.5000 \\
6&1.2000&1.4000&1.6000&1.8000&2.0000&2.2000&2.2333&2.2667&2.3000&2.3333 \\
7&1.1667&1.3333&1.5000&1.6667&1.8333&2.0000&2.1667&2.1905&2.2143&2.2381 \\
8&1.1429&1.2857&1.4286&1.5714&1.7143&1.8571&2.0000&2.1429&2.1607&2.1786 \\
9&1.1250&1.2500&1.3750&1.5000&1.6250&1.7500&1.8750&2.0000&2.1250&2.1389 \\
10&1.1111&1.2222&1.3333&1.4444&1.5556&1.6667&1.7778&1.8889&2.0000&2.1111 \\
11&1.1000&1.2000&1.3000&1.4000&1.5000&1.6000&1.7000&1.8000&1.9000&2.0000
\end{array}
$$}
\caption[]{\footnotesize
Values of $\kappa_S(n)$ for $S=\{1/B\cd 1/B)\}$ with $2\le B\le 11$
and $1\le n\le 10$.}
\label{kappa-table}
\end{table}

\begin{figure}[ht]
\centering\mbox{
\vbox{\psfig{file=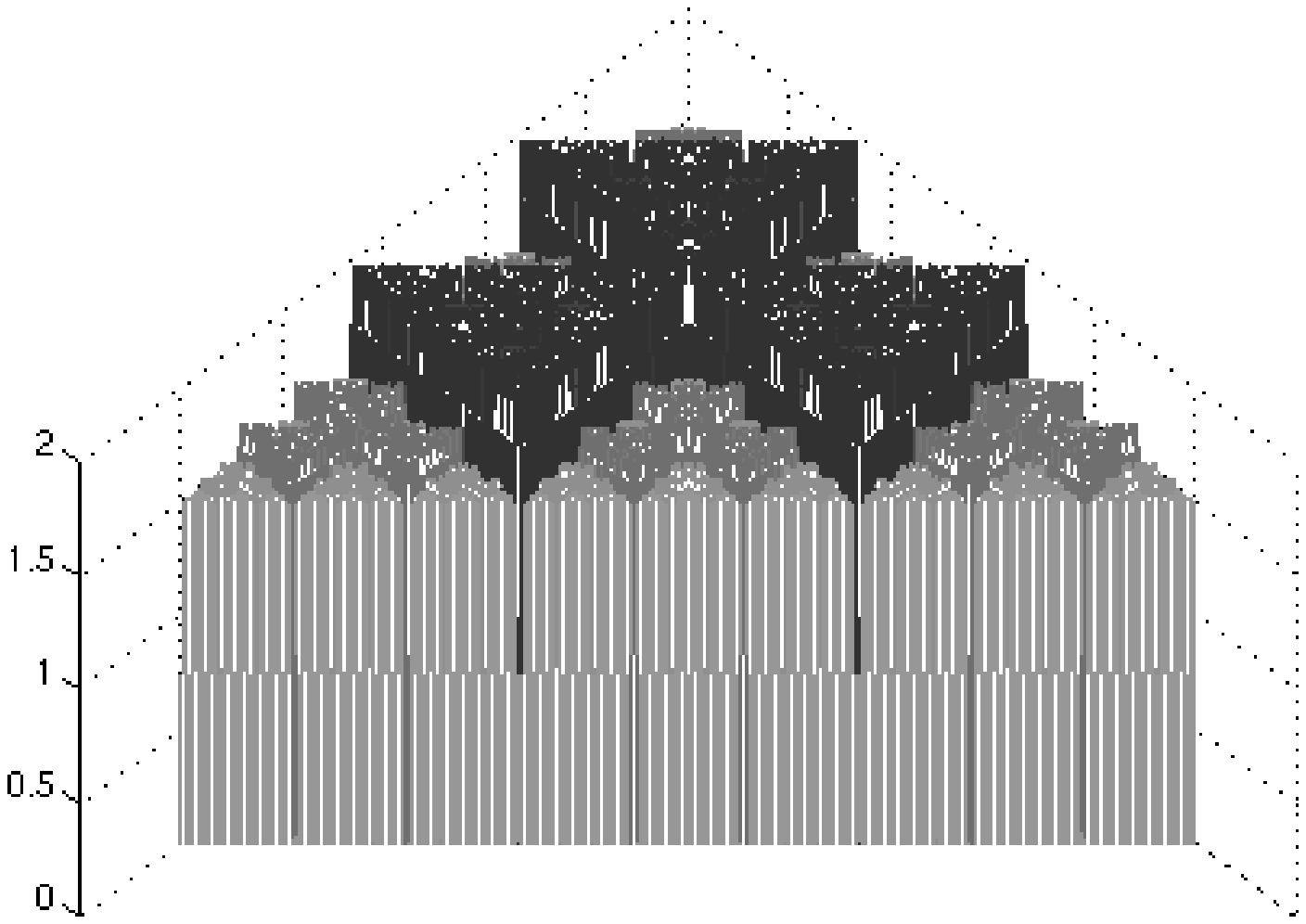,width=3.2in}\psfig{file=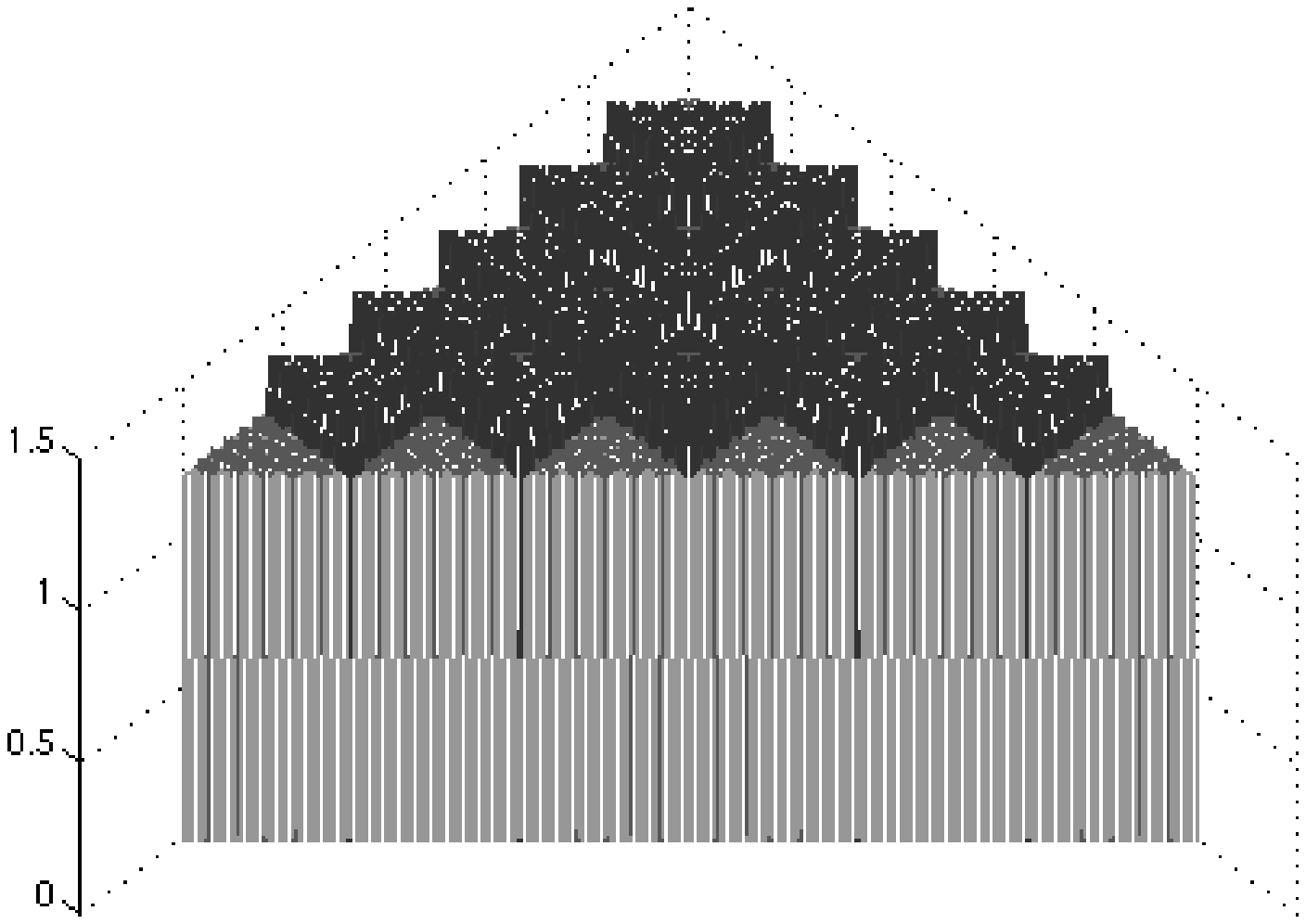,width=3.2in}\psfig{file=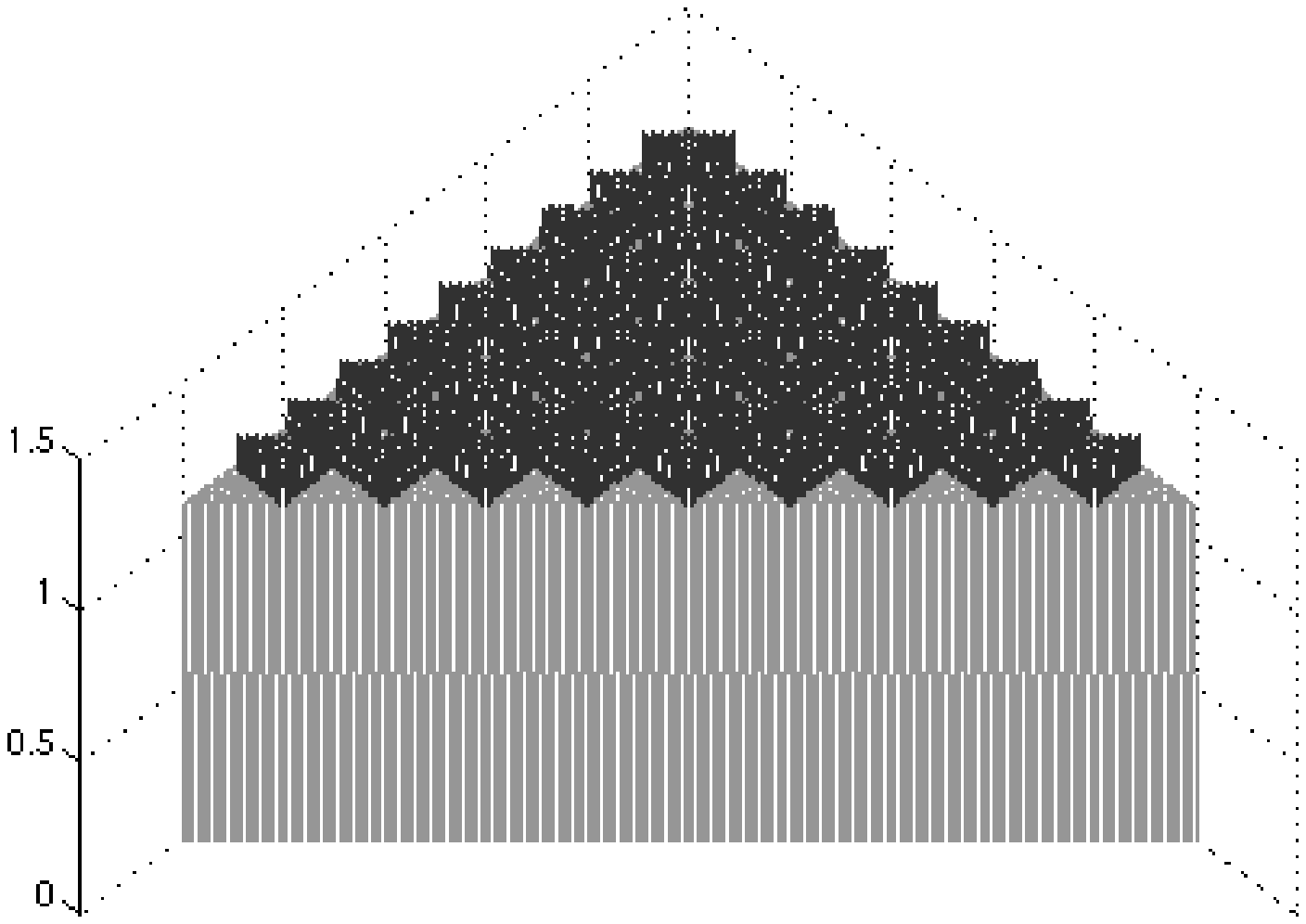,width=3.2in}}}
\caption[]{\footnotesize Graphs of $z=E_S(x,y,1-x-y)$ for $S=\{(1/B\cd
1/B)\}\in \Delta_{B-1}$ showing the dependence on $B$. The values of
$B$ are $B=3$ (top), $B=6$ (middle), and $B=10$ (bottom).  (The graph
for $B=2$ is in \cite[p.~23]{Dilworth-Howard-Roberts:1}.) }
\label{E-graph}
\end{figure}

The graphs of $z=E^{\Delta_2}_S(x,y,1-x-y)$ for some small values of
$B$ are given in Figure~\ref{E-graph}.

\begin{remark}\label{n=1}
Let $\mathcal{B}$ be the set of numbers of the form $j/B^l$ for $j,l$
integers and $j$ relatively prime to $B$.  Then using the series
expansion\eq{H-fp} and the argument of \cite[Prop.~2.25
p.~22]{Dilworth-Howard-Roberts:1} it is not hard to show
$E=E^{\Delta_1}_S\cn\Delta_1\to \R$ is given by
$$
E(1-t,t)=\begin{cases} \dfrac{B}{B-1},& t\notin\mathcal{B};\\
\dfrac{B}{B-1}-\dfrac{1}{B^{l-1}},& 
	t=\dfrac{j}{B^{l-1}}\in\mathcal{B}.\end{cases}
$$
\null\vskip-22pt\qed
\end{remark}

\subsection{The formula for $E^{\Delta_n}_S$.}

Let $[B]=\{1,2\cd B\}$ and let $X=[B]^\N$.  Let $\mu$ be the measure
on $X$ given by $\mu=\prod_{j=1}^\infty \nu_j$ where $\nu_j$ is the
measure on $[B]$ given by $\mu_j(\{i\})=1/B$ for $1\le i\le B$.
Therefore if $I\in \pi_k$ then $\mu(I)=1/B^k$.  The following lemma on
being able to realize certain sequences of numbers as sequences $\la
\mu(I_j)\ra_{j=1}^\infty$ with $\la I_j\ra_{j=1}^\infty$ a sequence
from $\{\nothing\}\cup\bigcup_{k=0}^\infty\pi_k$ allows us to simplify
the definition of $E_S(x)$ in some cases by replacing the infimum over
$S$-ranked measures with an infimum over special sequences of numbers
rather than measures.

\begin{lemma}\label{seq}
Let $\la r_j\ra_{j=1}^\infty$ be a nondecreasing sequence of nonnegative
integers such that 
$$
\sum_{j=1}^\infty \frac{1}{B^{r_j}}\le 1.
$$
Then there is a disjoint sequence $\la I_j\ra_{j=1}^\infty$ in
$\bigcup_{k=0}^\infty \pi_k$ such that
$$
\mu(I_j)=\frac{1}{B^{r_j}}\quad \text{and}\quad r(I_j)=r_j.
$$
\end{lemma}

\begin{proof}
Since for $I\in \bigcup_{k=0}^\infty\pi_k$ we have $\mu(I)=1/B^{r(I)}$ it is
enough to show the existence of a disjoint sequence $\la
\mu(I_j)\ra_{j=1}^\infty$ with $\mu(I_j)=1/B^{r_j}$ for then
$r(I_j)=r_j$ automatically holds.  We select this sequence
recursively.  Suppose that $I_1,I_2\cd I_j$ have been chosen to be
pointwise disjoint with $\mu(I_i)=1/B^{r_i}$.  Then
$$
\sum_{i=1}^j\mu(I_i)=\sum_{i=1}^j\frac{1}{B^{r_i}}
\le 1-\sum_{i=j+1}^\infty\frac{1}{B^{r_i}}\le 1-\frac{1}{B^{r_{j+1}}}.
$$
Since each of the sets $I_1,I_2\cd I_j$ is a union of atoms from
$\pi_{r_{j+1}}$, there is an atom of $\pi_{r_{j+1}}$ that is disjoint
from $I_1,I_2\cd I_j$. As atoms of $\pi_{r_{j+1}}$ have $\mu$-measure
$1/B^{r_{j+1}}$ we can use this atom as $I_{j+1}$.
\end{proof}

In light of Lemma~\ref{seq} and Proposition~\ref{prop4}
the value of $E=E^{\Delta_n}_S$ at $x=(x_0\cd
x_n)\in\Delta_n$ is given by
\begin{align*}
E(x)&
=\min\left\{\sum_{j=1}^\infty r(I_j)\mu(I_j) : 
	x\mid \la \mu(I_j)\ra_{j=1}^\infty\right\}\\
\intertext{(where $\mu$ is $S$-ranked, $\la I_j\ra_{j=1}^\infty$ is
pairwise disjoint, and $\sum_{j=1}^\infty \mu(I_j)=1$)}
&=\min \left\{\sum_{j=1}^\infty \frac{r_j}{B^{r_j}} : 
	x\mid \la 1/B^{r_j}\ra_{j=1}^\infty\right\}\\
\intertext{(where $r_j\in \N$ and $\sum_{j=1}^\infty 1/B^{r_j}=1$)}
&=\sum_{k=1}^n\min_{\begin{subarray}{c}N_0\cd N_n\\ \text{partitions
}\N\end{subarray}} \left\{\sum_{j=1}^\infty \frac{r_j}{B^{r_j}} : 
x_k= \sum_{j\in N_j} \frac{1}{B^{r_{j}}}\right\}.
\end{align*}
So if $H\cn [0,1]\to \R$ is defined by $H(0)=0$ and 
$$
H(x)=H_B(x)=\min\left\{\sum_{j=1}^\infty \frac{r_j}{B^{r_j}}:
\sum_{j=1}^\infty \frac{1}{B^{r_j}}=x\right\}
$$
for $x\in (0,1]$, then
$$
E(x)=\sum_{k=0}^nH(x_k).
$$
(We will shortly see that $H_B$ is also given by the formula\eq{H-fp}
so this notation is consistent with the notation used in the statement
of Theorem~\ref{main-bary}.)

We now give some other representations of $H$.  For
$x\in[0,1]$ consider sums
$$
\sum_{j=1}^\infty \frac{r_j}{B^{r_j}}\quad \text{where}\quad
x=\sum_{j=1}^\infty \frac{1}{B^{r_j}}.
$$ 
Let $x_i=|\{j: r_j=i\}|$.  Then these sums can be rewritten as
$$
\sum_{i=0}^\infty \frac{ix_i}{B^i}\quad \text{where}\quad
x=\sum_{i=0}^\infty \frac{x_i}{B^{i}}
$$
and so
\begin{equation}\label{B-adic-H}
H(x)=\min\left\{ \sum_{i=0}^\infty \frac{ix_i}{B^i} : \sum_{i=0}^\infty
\frac{x_i}{B^i} =x,\ x_i\in \N\right\}.
\end{equation}

\begin{lemma}\label{xi-B}
If $\sum_{i=0}^\infty ix_i/B^i$ is a minimizing sum in\eq{B-adic-H}
(so that $H(x)=\sum_{i=0}^\infty ix_i/B^i$),
then $x_i\in\{0,1\cd B-1\}$. 
\end{lemma}

\begin{proof}
Clearly $x_0\le 1$ (otherwise $x\notin[0,1]$).  Suppose that for some
$j\ge 1$ that $x_j\ge B$. Then let
$$
y_i=\begin{cases}x_{j-1}+1,& i=j-1;\\ x_j-B,& i=j;\\ x_i & i\ne j,j-1.
	\end{cases}
$$
Then each $y_i$ is nonnegative integer, $\sum_{i=0}^\infty y_i/B^i=x$
and 
\begin{align*}
\sum_{i=0}^\infty\frac{iy_i}{B^i}&=\frac{(j-1)(x_{j-1}+1)}{B^{j-1}}
+\frac{j(x_j-B)}{B^j}+\sum_{i\ne j,j-1}\frac{ix_i}{B^i}\\
&=\frac{j-1}{B^{j-1}}-\frac{j}{B^{j-1}}+\frac{(j-1)x_{j-1}}{B^{j-1}}
+\frac{jx_j}{B^j}+\sum_{i\ne j,j-1}\frac{ix_i}{B^i}\\
&=-\frac{1}{B^{j-1}}+\sum_{i=0}^\infty\frac{ix_i}{B^i}=H(x)-\frac{1}{B^{j-1}}. 
\end{align*}
This contradicts the minimality of the sum and completes the proof.
\end{proof}

Recall that any real number $x\in [0,1]$ has a base $B$-expansion
$x=\sum_{i=0}^\infty x_i/B^i$ where each $x_i\in\{0,1\cd B-1\}$.  This
expansion is unique unless $x$ is a $B$-adic rational (that is a
rational number of the form $k/B^l$ for integers $k$ and $l$).  A
$B$-adic rational has exactly two base $B$ expansions, one finite and
one infinite (if $x_n>0$ then
$\sum_{i=0}^nx_i/B^i=\sum_{i=0}^{n-1}x_i/B^i +(x_n-1)/B^n+
\sum_{i=n+1}(B-1)/B^i$).  For $B$-adic rationals $x$ we will always
use the finite expansion, but will still write $x=\sum_{i=0}^\infty x_i/B^i$
with the understanding that $x_i=0$ for $i$ sufficiently large.

\begin{prop}\label{base-B}
If $x\in [0,1]$ has base $B$ expansion 
$x=\sum_{i=0}^\infty x_i/B^i$, then $H(x)$ is given by
$$
H(x)=\sum_{i=0}^\infty \frac{ix_i}{B^i}.
$$
\end{prop}

\begin{proof}
From Lemma~\ref{xi-B} we know that if $x=\sum_{i=0}^\infty y_i/B^i$
with $y_i$ nonnegative integers is the expansion of $x$ so that $H(x)=
\sum_{i=0}^\infty iy_i/B^i$, then $0\le y_i\le B-1$.  When $x$ is not
a $B$-adic rational uniqueness of base $B$ expansions implies that
$y_i=x_i$ and we are done.  If $x$ is a $B$-adic rational and so has
two expansions with $0\le y_i\le B-1$ then direct
calculation shows that $\sum_{i=0}^\infty iy_i/B^i$
is smaller when the finite expansion is used.  Thus $y_i=x_i$ in this
case also.
\end{proof}

\begin{figure}[ht]
\centering
\mbox{\psfig{file=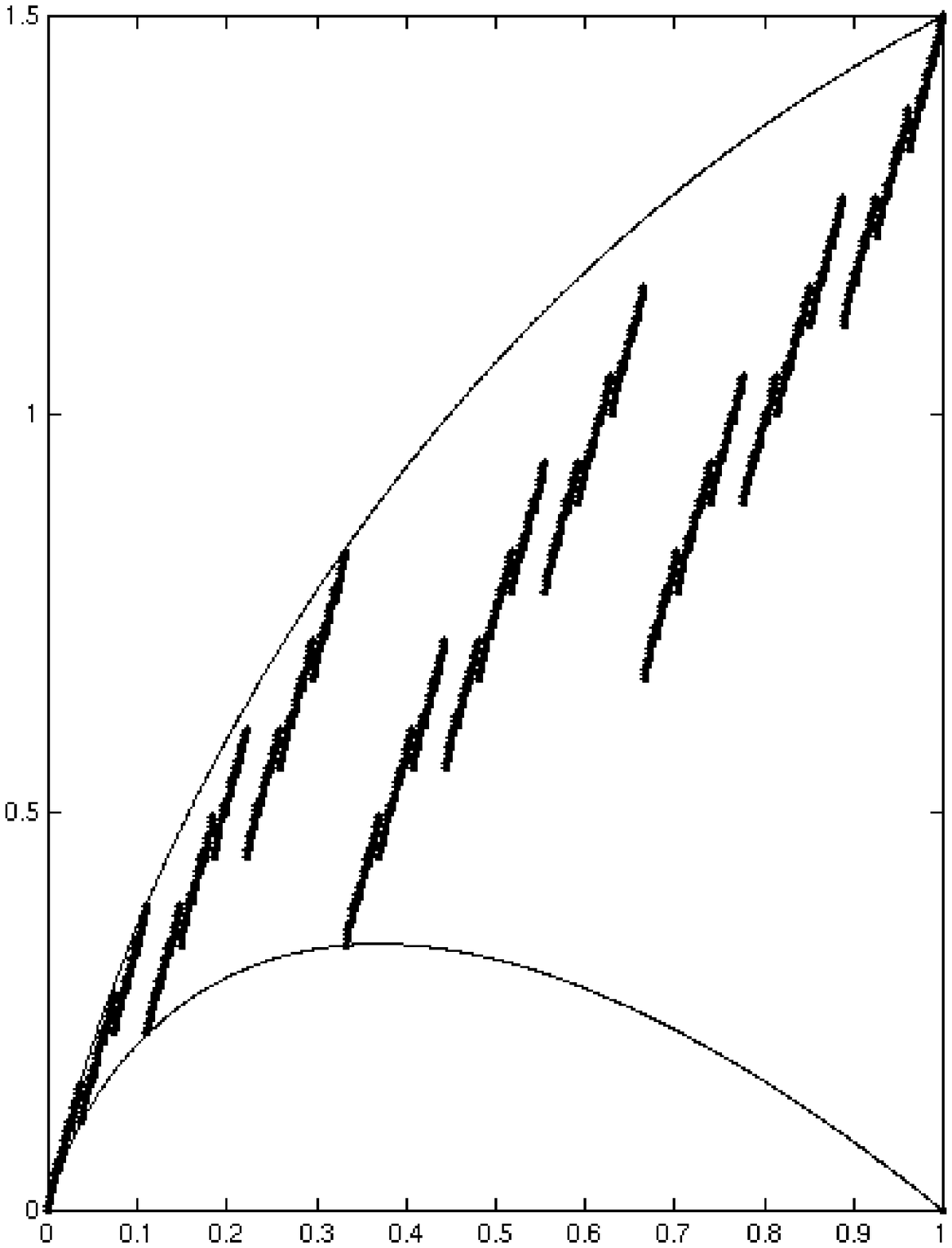,width=2.3in}\ \psfig{file=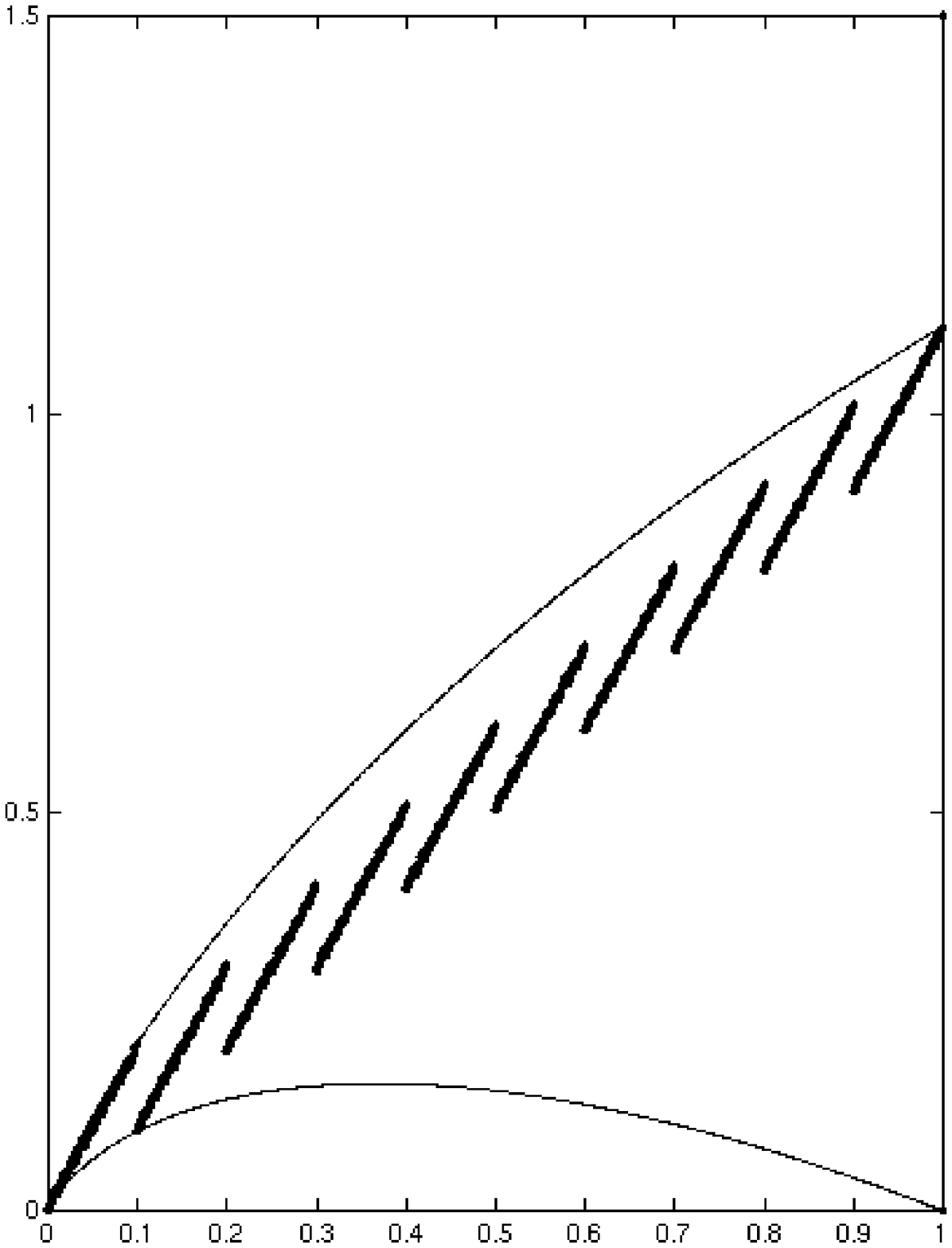,width=2.3in}}
\caption[]{\footnotesize
Graphs of $y=H_B(x)$, $y=x\log_B(1/x)$, and
$y=Bx/(B-1)+x\log_B(1/x)$ on $[0,1]$ for the values $B=3$ and $B=10$.}
\label{H-graph}
\end{figure}

It is convenient to extend $H$ to all of $\R$ to be periodic,
$H(x+1)=H(x)$.  This is possible as $H(0)=H(1)=0$.  Let $r\cn\R\to\R$
be the function that agrees with the greatest integer (or floor)
function on $[0,B)$ and is periodic of period~$B$. That is 
$$
r(x):=\begin{cases} \lfloor{x}\rfloor,& 0\le x< B;\\ r(x+B)=r(x),&
x\in\R.\end{cases}
$$
Then if $x=\sum_{i=1}^\infty x_i/B^i$ is the base $B$ expansion of
$x\in [0,1)$ then it is easily checked that $x_i=r(B^ix)$ and therefore
$x=\sum_{i=1}^\infty r(B^ix)/B^i$. Then the fractional part  $\fp{x}$ of the
real number $x$ is given by
$$
\fp{x}=\sum_{i=1}^\infty \frac{r(B^ix)}{B^i}
$$
as both sides are equal to $x$ on $[0,1)$ and are periodic of
period~$1$.  Also the periodic extension of $H$ to $\R$ is given by
$$
H(x)=\sum_{i=1}^\infty \frac{ir(B^ix)}{B^i}.
$$
These relations can be used to prove:
\begin{prop}\label{H-props}
The periodic extension of $H$ to $\R$ satisfies the functional
equation 
\begin{equation}\label{H-eqn}
H(x)=\fp{x} +\frac1 BH(Bx)
\end{equation}
and has the series representation
\begin{equation}\label{H-series}
H(x)=\sum_{k=0}^\infty \frac{\fp{B^kx}}{B^k}.
\end{equation}
Thus $H$ is lower semi-continuous, continuous at all points of
$[0,1]$ that are not $B$-adic rationals, and right continuous at all
points of $[0,1]$.  Also this function satisfies
the bounds
$$
x\log_B(1/x)\le H(x)\le \frac{Bx}{B-1}+x\log_B(1/x)
$$
on $[0,1]$ (see Figure~\ref{H-graph}).
\end{prop}

\begin{proof}
Other than the lower bound $x\log_B(1/x)\le H(x)$, we refer the reader
to the proofs of~\cite[Prop.~2.14 p.~15]{Dilworth-Howard-Roberts:1}
and~\cite[Prop.~2.21 p.~19]{Dilworth-Howard-Roberts:1} which cover the
case when $B=2$. Only trivial changes are required for the general
case.  

To prove the lower bound, suppose $x=\sum_{j=n}^\infty x_j/B^j$ is the
base $B$ expansion for $x$ with $x_n\ge1$.  Then $x\ge 1/B^n$ and
therefore $\log_B(1/x)\le n$.  Thus
$$
x\log_B(1/x)\le \sum_{j=n}^\infty\frac{nx_j}{B^j}\le \sum_{j=n}^\infty
\frac{jx_j}{B^j}=H(x)
$$
as required.
\end{proof}

We have now finished all of the proof of Theorem~\ref{main-bary} other
than computing the exact value of $\kappa_S(n)$.

\begin{remark}\label{H-self-cong}
The graph of $H_B$ has some interesting geometric properties.  The
following facts can be verified by the arguments used
in~\cite[Remark~2.15 p.~16]{Dilworth-Howard-Roberts:1} which
corresponds to the case $B=2$.  For all positive integers $i,j,k$ the
graphs of the restrictions $H_B\big|_{[i/B^k,(i+1)/B^k)}$ and
$H_B\big|_{[j/B^k,(j+1)/B^k)}$ are translates of each other and so the
graph of $H$ is ``locally self congruent at all scales $1/B^k$''.  The
closure of the graph is homeomorphic to the Cantor set and the graph
itself is this Cantor set with a countable number of points
deleted. Thus the graph is zero-dimensional as a topological space.
However the Hausdorff dimension of the graph is one.  Thus the closure
of the graph has metric dimension larger than its topological
dimension and therefore is a fractal.\qed
\end{remark}

\begin{remark}\label{E-self-sim}
(Cf.~\cite[Remark~2.26 p.~22]{Dilworth-Howard-Roberts:1})
The functional equation\eq{H-eqn} for $h=H_B$ can be used to explain the
self-similarities of the graph of $E_S\cn \Delta_n\to \R$ with
$S=\{(1/B\cd 1/B)$.  Let $v\in \Delta_{n}$ be a point so that all the
entries of $(B-1)v$ are integers.  Let $x\in \Delta_n$ be any point
that is not a vertex.  Then $(x+(B-1)v)/B$ is not a vertex and so all
the components of $(x+(B-1)v)/B$ are in the interval $[0,1)$ and thus
are equal to their fractional part.  So letting $x=(x_0\cd x_n)$ and
$v=(v_0\cd v_n)$ and using\eq{H-series}
\begin{align*}
E&\(\frac{x+(B-1)v}{B}\)
=\sum_{k=0}^nH\(\frac{x_k+(B-1)v_k}{B}\)\\
&=\sum_{k=0}^n\pmb{\bigg{\{}}\frac{x_k+(B-1)v_k}{B}\pmb{\bigg\}}
+\frac1B\sum_{k=0}^nH(x_k+(B-1)v_k)\\
&=\sum_{k=0}^n\frac{x_k+(B-1)v_k}{B}+\frac1B\sum_{k=0}^nH(x_k)\\
&=1+\frac1BE(x).
\end{align*}
where we have used the fact that for each $k$ such that $(B-1)v_k$ is an
integer that $H(x_k+(B-1)v_k)=H(x_k)$ as $H$ has period one.  On
the set $\Delta_n\times [0,\infty)$, for each $v\in \Delta_n$ such that
$(B-1)v$ has all integer entries, define $\theta_v\cn \Delta_n\times
[0,\infty)\to \Delta_n\times [0,\infty)$ by
$$
\theta_v(x,z)=\(\frac{x+(B-1)v}{B},1+\frac1B z\).
$$
This is the dilation by $1/B$ with center $(v,B/(B-1))$.  The
calculation we have just done shows for each $x\in \Delta_n$ that is
not a vertex that
\begin{align*}
\theta_v(x,E(x))&=\(\frac{x+(B-1)v}{B},1+\frac1BE(x)\)\\
&=\(\frac{x+(B-1)v}{B},E\(\frac{x+(B-1)v}{B}\)\).
\end{align*}
Therefore each of these dilations maps the graph of $E$ into a subset
of the graph.  When $B$ is much larger than $n$ there will be a
large number of points $v\in \Delta_n$ so that $(B-1)v$ has all
integral elements and thus in this case the graph of $z=E(x)$ will have
a very large number of self symmetries.  This is apparent in the
bottom graph in Figure~\ref{E-graph}  where $n=2$ and $B=10$.\qed
\end{remark}

\subsection{Calculation of $\kappa_S(n)$.}\label{sec:kappa-bary}
Let $\mathcal{B}_n$ be the points in $\Delta_n$ with $B$-adic rational
coordinates.  Then $\mathcal{B}_n$ is dense in $\Delta_n$ and $E$ is
lower semi-continuous.  Therefore
$$
\sup_{x\in\mathcal{B}_n}E(x)=\sup_{x\in\Delta_n}E(x).
$$
So there is a sequence $\la x(s)\ra_{s=1}^\infty\subset\Delta_n$ so
that $x(s)=\sum_{k=0}^nx_k(s)e_k$ with each $x_k(s)$ a $B$-adic
rational and with $\lim_{s\to\infty}E(x(s))=\kappa_S(n)$.  Each
$x_k(s)$ can be written $x_k(s)=\sum_{j=0}^\infty x_{j\,k}(s)/B^j$
with $x_{j\,k}(s)\in\{0\cd B-1\}$ and each sequence $\la
x_{j\,k}(s)\ra_{j=0}^\infty$ eventually~$0$.  By passing to a
subsequence we may assume that for $0\le k\le n$ and $0\le j <\infty$
that $\lim_{s\to \infty}x_{j\,k}(s)=x_{j\,k}$ with $x_{j\,k}\in\{0\cd
B-1\}$.  That is for fixed $j$ and $k$ we have $x_{j\,k}(s)=x_{j\,k}$
for sufficiently large $s$.  Therefore if $x_k=\sum_{j=0}^\infty x_{j\,
k}/B^j$ for $0\le k\le n$, then by the Lebesgue Dominated Convergence
Theorem $\sum_{k=0}^nx_k=1$.  (All the series $\sum_{j=0}^\infty
x_{j\,k}(s)/B^j$ are dominated by the convergent geometric series
$\sum_{j=0}^\infty(B-1)/{B^j}$ so we can take the limit, i.e.,
$1=\lim_{s\to\infty}\sum_{k=0}^nx_{j\,k}(s)/B^j
=\sum_{k=0}^nx_{j\,k}/B^j=\sum_{k=0}^nx_k$.)  Another application of
the Lebesgue Dominated Convergence
Theorem gives
$$
\kappa_S(n)=\lim_{s\to\infty}E(x(s))=\lim_{s\to\infty}\sum_{k=0}^\infty
\sum_{j=0}^\infty\frac{jx_{j\,k}(s)}{B_j}=\sum_{k=0}^\infty
\sum_{j=0}^\infty\frac{jx_{j\,k}}{B_j}.
$$
Let
$$
M_j(s):=\sum_{k=0}^nx_{j\,k}(s)\quad\text{and}\quad M_j:=\sum_{k=0}^nx_{j\,k}.
$$
So for fixed $j$ we have $M_j(s)=M_j$ for sufficiently large $s$.
Also
$$
E(x(s))=\sum_{j=0}^\infty \frac{jM_j(s)}{B^j}, \quad
\kappa_S(n)=\sum_{j=0}^\infty \frac{jM_j}{B^j},\quad 
1=\sum_{j=1}\frac{M_j}{B^j},
$$
and for fixed $s$ we have $M_j(s)=0$ for sufficiently large
$j$.

As a first observation note that each $x_j(s)\le B-1$ which implies
$M_j(s)\le (n+1)(B-1)$ which in turn implies
\begin{equation}\label{(1)}
M_j\le (n+1)(B-1).
\end{equation}

Assuming $n\ge 1$ (obviously $\kappa_S(0)=0$) we have $M_0=0$ (for 
$M_0=\sum_{k=0}^nx_{0\,k}>0$ would imply that the point $(x_0\cd x_n)$
is a vertex of $\Delta_n$ and this is clearly not a maximizing
sum).
Let
$$
\ell+1=\text{least $j$ such that $M_j>0$.}
$$
In particular $0=M_0=\cdots =M_\ell$ and $M_{\ell+1}>0$.

\begin{lemma}\label{(2)}
If $j\ge \ell+2$, then $M_j\ge (B-1)n$.
\end{lemma}

\begin{proof}
Suppose not and let $i$ be the least $i\ge \ell+2$ such that
$M_i<(B-1)n$.  If $i>\ell+2$, then $M_{i-1}\ge (B-1)n$ and if
$i=\ell+2$, then $M_{i-1}=M_{\ell+1}>0$.  In either case $M_{i-1}>0$.
There is an $s_0$ such that for $s\ge s_0$, $M_{i-1}(s)=M_{i-1}$ and
$M_i(s)=M_i$.  Thus for each $s\ge s_0$ there is a
$y(s)=\sum_{k=0}^n\(\sum_{j=0}^\infty y_{j\,k}(s)/B^j\)e_k$ with
$y_{j\,k}(s)$ defined so that
$$
y_{j\,k}(s)=x_{j\,k}(s)\quad\text{if}\quad j\ne i-1,i,
$$
$$
\sum_{k=0}^ny_{i-1\, k}(s)=M_{i-1}-1.
$$
(this is possible because $M_{i-1}>0$) and
$$
\sum_{k=0}^ny_{i\,k}(s)=M_i+B
$$
(this is possible becasue $M_i<(B-1)n$ so that $M_i+B\le (B-1)(n+1)$).
But then for $s>s_0$,
$$
E(y(s))=E(x(s))+\frac{iB}{B^i}-\frac{i-1}{B^{i-1}}=E(x(s))+\frac{1}{B^{i-1}}.
$$
But then $\lim_{s\to\infty}E(y(s))=\kappa_S(n)+1/B^{i-1}$ which is impossible.
\end{proof}

\begin{lemma}\label{(3)} For infinitely many $j$ the inequaltiy
$M_j<(B-1)(n+1)$ holds.
\end{lemma}

\begin{proof}
Suppose that for some $j_0$ that $j\ge j_0$ implies $M_j=(B-1)(n+1)$.
Then there exists $s_0$ such that for $j<j_0$ and $s>s_0$ we have
$M_j(s)=M_j$.  But then for any $s>s_0$ (recall that for fixed $s$
there holds $M_j(s)=0$ for $j$ sufficiently large)
\begin{align*}
1&=\sum_{j=0}^\infty\frac{M_j(s)}{B^j}=\sum_{j=0}^{j_0-1}\frac{M_j(s)}{B^j}
	+\sum_{j=j_0}^\infty\frac{M_j(s)}{B^j}\\
&=\sum_{j=0}^{j_0}\frac{M_j}{B^j}+\sum_{j=j_0}^\infty\frac{M_j(s)}{B^j}
<\sum_{j=0}^{j_0}\frac{M_j}{B^j}+\sum_{j=j_0}^\infty\frac{(B-1)(n+1)}{B^j}\\
&=\sum_{j=0}^\infty \frac{M_j}{B^j}=1
\end{align*}
which is a contradiction.
\end{proof}

\begin{lemma}\label{(4)}
If $j\ge \ell+2$, then $M_j=(B-1)n$.
\end{lemma}

\begin{proof}
By Lemma~\ref{(2)} $M_j\ge (B-1)n$ and by Lemma~\ref{(3)},
$M_j<(B-1)(n+1)$ for infinitely many $j$.  Thus
\begin{align*}
\sum_{j=\ell+2}^\infty\frac{M_j}{B^j}&=\sum_{j=\ell+2}^\infty\frac{(B-1)n}{B^j}
	+\sum_{j=\ell+m+2}^\infty\frac{M_j-(B-1)n}{B^j}\\
&=\frac{(B-1)n}{B^{\ell+2}}\(\frac{1}{1-1/B}\)
	+\sum_{j=\ell+2}^\infty\frac{M_j-(B-1)n}{B^j}\\
&=\frac{n}{B^{\ell+1}}+\sum_{j=\ell+2}^\infty\frac{M_j-(B-1)n}{B^j}
\end{align*}
Set $R=\sum_{j=\ell+2}^\infty\frac{M_j-(B-1)n}{B^j}$.  Then
$$
0\le R< \sum_{j=\ell+2}^\infty
\frac{B-1}{B^j}=\frac{B-1}{B^{\ell+2}}\(\frac{1}{1-1/B}\) =\frac{1}{B^{\ell+1}}.
$$
where the first inequality follows form Lemma~\ref{(2)} and the second
from Lemma~\ref{(3)}.  Thus
$$
1=\sum_{j=0}^\infty\frac{M_j}{B^j}=\sum_{j=0}^{\ell+1}\frac{M_j}{B^j}+\frac{n}{B^{\ell+1}}+R
=\frac{L}{B^{\ell+1}}+R
$$
with $0\le R<1/B^{\ell+1}$ and $L$ a positive integer.  But then $0\le
1- L/B^{\ell+1} =R<1/B^{\ell+1}$ which implies $R=0$.  That is
$0=R=\sum_{j=\ell+2}^\infty \frac{M_j-(B-1)n}{B^j}$.  Thus
$M_j-(B-1)n=0$ for $j\ge\ell+2$.
\end{proof}

\begin{lemma}\label{(5)}
The integer $\ell$ satisfies $\dfrac{M_{\ell+1}+n}{B^{\ell+1}}=1$.
\end{lemma}

\begin{proof}  Using the results from the last several lemmas:
\begin{align*}
1&=\sum_{j=0}^\infty\frac{M_j}{B^j}=\frac{M_{\ell+1}}{B^{\ell+1}}+(B-1)n\sum_{j=\ell+2}^\infty\frac{1}{B^j}\\
&=\frac{M_{\ell+1}}{B^{\ell+1}}+\frac{n}{B^{\ell+1}}
=\frac{M_{\ell+1}+n}{B^{\ell+1}}.
\end{align*}
\null\vskip-22pt
\end{proof}

\begin{lemma}\label{(6)}
The integer $\ell$ satisfies $B^\ell\le n< B^{\ell+1}$ so that
$\ell=\lfloor\log_Bn\rfloor$.
\end{lemma}

\begin{proof}
By Lemma~\ref{(5)} $M_{\ell+1}+n=B^{\ell+1}$ and $M_{\ell+1}>0$ so
$n<B^{\ell+1}$. For the other inequality, use
$M_{\ell+1}\le(n+1)(B-1)$ so that
\begin{align*}
B^{\ell+1}&=M_{\ell+1}+n\le (n+1)(B-1)+n=nB+B-1\\
\implies&\quad(n+1)B\ge B^{\ell+1}+1\\
\implies &\quad(n+1)B>B^{\ell+1}\\
\implies &\quad n+1> B^\ell\\
\implies &\quad n\ge B^\ell
\end{align*}
\null\vskip-22pt
\end{proof}

Using the results of these lemmas we can now compute the value of
$\kappa_S(n)$.  
\begin{equation}\label{pre-kappa}
\kappa_S(n)=\sum_{j=0}^\infty\frac{jM_j}{B^j} 
	=\frac{(\ell+1)M_{\ell+1}}{B^{\ell+1}}+
	n(B-1)\sum_{j=\ell+2}^\infty\frac{j}{B^j}.
\end{equation}
Using Lemma~\ref{lemma:sum} (with $x=1/B$ and $a=k=\ell+2$)
\begin{align*}
(B-1)\sum_{j=\ell+2}^\infty\frac{j}{B^j}&=(B-1)\sum_{i=0}^\infty\frac{\ell+2+i}{B^{\ell+2+i}}\\
&=(B-1)\frac{(\ell+2)(1/B)^{\ell+2}+(1-(\ell+2))(1/B)^{\ell+3}}{(1-1/B)^2}\\
&=\frac{(\ell+2)B-(\ell+1)}{(B-1)B^{\ell+1}}
\end{align*}
Substituting this and also $M_{\ell+1}=B^{\ell+1}-n$ (Lemma~\ref{(5)})
into\eq{pre-kappa} gives
\begin{align*}
\kappa_S(n)&=\frac{(\ell+1)M_{\ell+1}}{B^{\ell+1}}
	+\frac{n[(\ell+2)B-(\ell+1)]}{(B-1)B^{\ell+1}}\\
&=\frac{(\ell+1)(B^{\ell+1}-n)}{B^{\ell+1}}
	+\frac{n[(\ell+2)B-(\ell+1)]}{(B-1)B^{\ell+1}}\\
&=\ell+1+ \frac{-n(\ell+1)}{B^{\ell+1}}
	+\frac{n[(\ell+2)B-(\ell+1)]}{(B-1)B^{\ell+1}}\\
&=\ell+1+\frac{n}{(B-1)B^\ell}\\
&=\lfloor\log_B n\rfloor+1+\frac{n}{(B-1)B^{\lfloor\log_Bn\rfloor}}.
\end{align*}
This completes the proof of Theorem~\ref{main-bary}.


\providecommand{\bysame}{\leavevmode\hbox to3em{\hrulefill}\thinspace}
\providecommand{\MR}{\relax\ifhmode\unskip\space\fi MR }
\providecommand{\MRhref}[2]{%
  \href{http://www.ams.org/mathscinet-getitem?mr=#1}{#2}
}
\providecommand{\href}[2]{#2}

\end{document}